\documentclass[preprint,10pt]{elsarticle}
\usepackage{graphicx}
\usepackage[usenames, dvipsnames]{color}
\usepackage{caption}
\usepackage{amsfonts}
\usepackage{amsmath}
\usepackage[pdftex,dvipsnames]{xcolor}  
\usepackage{amsthm}
\usepackage{xargs}
\usepackage[colorinlistoftodos,prependcaption,textsize=tiny]{todonotes}
\newcommandx{\question}[2][1=]{\todo[linecolor=Plum,backgroundcolor=Plum!25,bordercolor=Plum,#1]{#2}}
\usepackage{makecell}
\usepackage{verbatim}
\usepackage{epsfig, algorithm}
\def\beq{\begin{equation}}
\def\eeq{\end{equation}}
\def\esplit{\end{split}}
\def\beqalign{\begin{array}{rl}}
\def\eeqalign{\end{array}}
\usepackage{amsfonts}

\usepackage{algorithm,algcompatible,amsmath}
\algnewcommand\INPUT{\item[\textbf{Input:}]}%
\algnewcommand\OUTPUT{\item[\textbf{Output:}]}%

\def\Rbb{\mathbb{R}}

\def\Abold{\mathbf{A}}
\def\Bbold{\mathbf{B}}
\def\Cbold{\mathbf{C}}
\def\Dbold{\mathbf{D}}

\def\Gbold{\mathbf{G}}
\def\Hbold{\mathbf{H}}
\def\Ibold{\mathbf{I}}

\def\Kbold{\mathbf{K}}

\def\Mbold{\mathbf{M}}
\def\Nbold{\mathbf{N}}

\def\Pbold{\mathbf{P}}
\def\Qbold{\mathbf{Q}}
\def\Rbold{\mathbf{R}}

\def\Ubold{\mathbf{U}}
\def\Vbold{\mathbf{V}}
\def\Wbold{\mathbf{W}}
\def\Xbold{\mathbf{X}}
\def\Ybold{\mathbf{Y}}
\def\Zbold{\mathbf{Z}}

\def\bbold{\mathbf{b}}
\def\cbold{\mathbf{c}}

\def\qbold{\mathbf{q}}
\def\pbold{\mathbf{p}}

\def\ubold{\mathbf{u}}

\def\wbold{\mathbf{w}}
\def\xbold{\mathbf{x}}

\def\lambdabold{\boldsymbol{\lambda}}

\def\mubold{\boldsymbol{\mu}}

\def\Sigmabold{\mathbf{\Sigma}}
\def\sigmabold{\mathbf{\sigma}}
\def\Gammabold{\mathbf{\Gamma}}
\def\Omegabold{\mathbf{\Omega}}

\def\zetabold{\boldsymbol{\zeta}}
\def\zerobold{{\bf 0}}
\def\0{\mathbf{\0}}

\newenvironment{@abssec}[1]{%
     \if@twocolumn
       \section*{#1}%
     \else
       \vspace{.05in}\footnotesize
       \parindent .2in
         {\upshape\bfseries #1. }\ignorespaces
     \fi}
     {\if@twocolumn\else\par\vspace{.1in}\fi}

\newcommand\keywordsname{Key words}
\newcommand{\bmat}[1]{\begin{bmatrix}#1\end{bmatrix}} 

\newtheorem{mydef}{Definition}

\usepackage{graphicx}
\usepackage{amssymb}
\usepackage{amsmath}
\usepackage{lineno}

\journal{Journal Name}

\begin{document}

\begin{frontmatter}

\title{Gradient-based Constrained Optimization Using a Database of Linear Reduced-Order Models}

\author[label1]{Youngsoo Choi}
\author[label2]{Gabriele Boncoraglio}
\author[label2]{Spenser Anderson}
\author[label5]{David Amsallem}
\author[label2,label3,label4]{Charbel Farhat}

\address[label1]{Computational Engineering Division, Lawrence Livermore National Laboratory\footnote{Lawrence Livermore National Laboratory is operated by Lawrence Livermore National Security, LLC, for the U.S. Department of Energy, 
National Nuclear Security Administration under Contract DE-AC52-07NA27344 and LLNL-JRNL-796678}, Livermore, CA 94550}
\address[label2]{Department of Aeronautics and Astronautics, Stanford University, Stanford, CA 94305}
\address[label3]{Department of Mechanical Engineering, Stanford University, Stanford, CA 94305}
\address[label4]{Institute for Computational and Mathematical Engineering, Stanford University, Stanford, CA 94305}
\address[label5]{Facebook, Menlo Park, CA 94025}

\begin{abstract}
A methodology grounded in model reduction is presented for accelerating the gradient-based solution of a family of linear or nonlinear constrained optimization problems where the constraints
include at least one linear Partial Differential Equation (PDE). 
A key component of this methodology is the construction, during an offline phase, of a database of pointwise, linear, Projection-based Reduced-Order 
Models (PROM)s associated with a design parameter space and the linear PDE(s). A parameter sampling procedure based on an appropriate saturation assumption is proposed to maximize the efficiency of such
a database of PROMs. A real-time method is also presented for interpolating at any queried but unsampled parameter vector in the design parameter space the relevant sensitivities of a PROM. The practical 
feasibility, computational advantages, and performance of the proposed methodology are demonstrated for several realistic, nonlinear, aerodynamic shape optimization problems governed by linear aeroelastic constraints.
\end{abstract}


\end{frontmatter}


\section{Introduction}
\label{sec:Introduction}

Partial Differential Equation (PDE)-based constrained optimization problems arise in numerous engineering applications including automatic design, optimal control, 
and inverse identification. In practice, such problems are typically solved using a Nested ANalysis and Design (NAND) approach, which incurs the solution of many 
instances of the constraining parametric PDE, for different sampled parameter vectors in the design parameter space denoted here
by $\mathcal D$. In many industrial settings, such a PDE may model, 
for example, a Computational Fluid Dynamics (CFD) or Computational Structural Dynamics (CSD) problem, its discretization is high-dimensional and therefore computationally intensive, 
and hence its enforcement can represent 
a significant fraction of the total cost of the optimization problem. Of course, this situation is exacerbated when the optimization problem of interest is constrained by multiple such PDEs.
For this reason, PDE-constrained optimization problems are often solved using a less computationally intensive surrogate model such as, for example, a set of response surfaces for some Quantities of Interest
(QoIs), or a Projection-based Reduced-Order Model (PROM). 

In general, a PROM seeks an approximate solution of the underlying High-Dimensional Model (HDM) within a subspace of the high-dimensional state space that is spanned by a Reduced-Order Basis (ROB)
\cite{moore1981principal, sirovich1987turbulence}. In the context of an optimization problem, an HDM is typically a 
$\mubold$-parametric HDM, where $\mubold \in \mathcal D \subset \Rbb^{N_{\mathcal D}}$ is a design parameter vector, and 
$N_{\mathcal D}$ denotes the dimension of the design parameter space $\mathcal D$. It is well-known that a PROM constructed for a 
specific design parameter vector 
$\mubold \in \mathcal D$, which is referred to in this paper as a {\it pointwise} PROM, 
does not necessarily perform well at another design parameter vector $\mubold^{\prime}$ in this space (for example, see 
\cite{epureanu2003parametric, lieu2004parameter}). To adress this issue, at least three approaches have been proposed so far. 

In the first approach, a global ROB is constructed offline such that its associated PROM is accurate in the entire design
parameter space, or at least a large region of this space 
\cite{leGresley2000airfoil, weickum2009multi, amsallem2013design}. A drawback of this approach, particularly for high-dimensional 
design parameter spaces, is that it may lead to a ROB that is 
prohibitively large, or sufficiently large to offset the computational advantages of the associated PROM. For nonlinear problems -- in this context, the discretized nonlinear PDE constraints -- the concept 
of a cluster of local, low-dimensional ROBs proposed in \cite{amsallem2012local} addresses this issue. In this concept, locality refers to a region of the solution manifold, and the most appropriate local ROB 
is selected to approximate the solution when it evolves on its manifold.

Unlike the first approach which is comprehensive, the second approach was motivated primarily by optimization problems. In this approach, an initial ROB is {\it progressively} adapted so that it remains
accurate not necessarily in the entire design parameter space, but in those regions of this space that are visited by the optimization algorithm 
\cite{fahl2003reduced, yue2013accelerating, zahr2013construction, zahr2014progressive}. Specifically, as new parameter vectors
are queried by the optimizer in the design parameter space, 
additional solution snapshots are computed
only when needed to maintain accuracy, and the dimension of the PROM is adjusted to incorporate the information gathered from these snapshots. This approach is suitable when the traditional offline-online
computational framework is not particularly advantageous, and the minimization of the overall computational cost is instead preferred. On the other hand, it is not necessarily advantageous when
the offline-online approach is desired, and the objective is to accelerate the execution of the online phase.

For linear problems -- and in the context of this paper, specifically for the discretized linear PDE constraints -- the parameter dependence of a PROM associated with
a $\mubold$-parametric HDM can be addressed using the concept of a database of PROMs introduced in \cite{amsallem2009method, amsallem2011online}. This approach can be summarized as follows.
First, the design parameter space of interest is sampled offline using an appropriate greedy procedure, and a linear, pointwise PROM is constructed at each sampled design parameter vector using any preferred model
reduction technique. Then, at each queried but unsampled design parameter vector, a PROM is constructed by interpolating the pre-computed PROMs on one or multiple matrix manifolds that characterize them. 
For linear problems, this approach is attractive for two reasons: the dimension of each pointwise PROM can be kept small while maintaining local accuracy, due to the restriction of the scope of such a PROM
to a single parameter vector in the design parameter space (and perhaps a relatively small neighborhood); and the aforementioned interpolation on matrix manifolds can be performed in real time. 

This paper focuses on the family of constrained optimization problems where the objective function may be linear or nonlinear but is time-independent, and the PDE constraint -- and in the case of 
multiple PDE constraints, at least on of them -- is linear and its discretization is high-dimensional. For such problems, which abound in engineering applications, the concept of a database of pointwise, linear
PROMs is particularly attractive when multiple optimization instances must be carried out as in robust optimization, multi-objective optimization, and global optimization with multiple initial solutions. 
The work described here advances the state of the art of this concept, particularly for the solution of the aforementioned class of PDE-constrained optimization problems by a gradient-based NAND approach. 
Specifically, it contributes a parameter sampling procedure based on a saturation assumption for the construction of a cost-efficient PROM database. It also contributes a novel algorithm for interpolating the 
sensitivities of a pointwise, linear PROM with respect to parameter variations. Finally, it demonstrates the potential of the concept of a database of pointwise, linear PROMs equipped with these contributions 
for the identified class of applications. For this last purpose, it reports on the application of the proposed computational framework to the aerodynamic shape optimization of a nonlinear aeroelastic wing under 
a linear stability (flutter) constraint -- a problem which is normally untractable when linearized CFD and CSD models are preferred for representing the PDE constraint.

\section{PDE-Constrained Optimization}
\label{sec:PDE-ConstrainedOptimization}

\subsection{Problem formulation}

Consider a time-independent but otherwise general PDE-constrained optimization problem, where the objective function may be linear or nonlinear and at least one of the parametric PDE constraints is linear. This 
problem may be also governed by additional linear or nonlinear algebraic constraints. The main objective of this paper is to present a computational framework for accelerating, in this context, the enforcement 
of any parametric, linear, PDE constraint whose discretization is high-dimensional. For simplicity, but without any loss of generality, assume that the optimization problem contains only one such PDE constraint 
(the extension of the proposed computational framework to multiple parametric, linear, PDE constraints is straightforward).

The general form of a discretized, time-independent, $\mubold$-parametric, linear PDE can be written as
\begin{equation}\label{eq:LinearForm}
	\Rbold_L \left(\wbold_L(\mubold), \mubold \right) = \Abold (\mubold) \wbold_{L} (\mubold) - \bbold (\mubold) = \mathbf{0}
\end{equation}
where: the design parameter vector $\mubold \in \mathcal D \subset \Rbb^{N_{\mathcal D}}$ is referred to throughout the remainder of this
paper as the parameter ``point'' in $\mathcal D$,
$\displaystyle{\Abold(\mubold) \in \Rbb^{N_{w_{L}} \times N_{w_{L}}}}$ and $\bbold(\mubold) \in \Rbb^{N_{w_{L}}}$ are the governing parametric matrix and right-hand side arising from the discretization
of this PDE by $ N_{w_{L}}$ degrees of freedom (dof)s; and $\wbold_L(\mubold) \in \Rbb^{N_{w_{L}}}$ denotes its parametric solution. The parameter point $\mubold$ may represent, for example, the material properties, shape design parameters, or boundary 
conditions of a physical system of interest. 

For the sake of simplicity, it is also assumed here and in the remainder of this paper that when $\mubold$ is varied, the size and topology of the computational mesh underlying the discrete equation 
(\ref{eq:LinearForm}) remain unchanged. The extension of the computational framework presented in this paper to the case where this size and topology vary when $\mubold$ is varied can be performed
using the computational approaches presented in \cite{amsallem2016realtime}.

The general PDE-constrained optimization problem outlined above can be formulated as follows
\begin{equation}\label{eq:generalPb}
      \boxed{
        \begin{aligned}
& \underset{\wbold_L, \wbold_{NL}, \mubold}{\text{min}}
& & f\left(\wbold_{L}(\mubold),\wbold_{NL}(\mubold),\mubold\right) \\
& \text{s.t.} & &  \cbold\left(\wbold_{L}(\mubold),\wbold_{NL}(\mubold),\mubold)\right) \leq \boldsymbol{0} \\
		& & &  \Rbold_L \left(\wbold_L(\mubold), \mubold \right) = \Abold(\mubold) \wbold_{L}(\mubold) - \bbold(\mubold) = \mathbf{0} \\
& & &  \Rbold_{NL}\left(\wbold_{NL}(\mubold),\mubold\right)= \mathbf{0} \\
 \end{aligned}
      }
    \end{equation}
where: $f$ denotes the objective function and is assumed here for simplicity to be a scalar function; $\wbold_{NL}  \in \Rbb^{N_{w_{NL}}}$ is the $N_{w_{NL}}$-dimensional solution vector
of the set of discretized nonlinear PDEs $\Rbold_{NL}(\cdot,\cdot)= \mathbf{0}$ that may be present or absent from (\ref{eq:generalPb}); and $\cbold(\cdot, \cdot)$ is a general set of linear and/or nonlinear 
algebraic constraints including, for instance, linear constraints on $\wbold_L$ and/or $\wbold_{NL}$ or box constraints on $\mubold$.
  
Because this paper focuses on the NAND solution approach, problem \eqref{eq:generalPb} is viewed here as an optimization problem defined over the parameter point $\mubold \in \mathcal D$ only; the state vectors 
$\wbold_L$ and $\wbold_{NL}$ are considered as functions of $\mubold$ that are implicitly defined by the discretized PDE constraints. When the optimizer queries a value of $\mubold$, an appropriate PDE solver is 
applied to compute each of $\wbold_L$ and $\wbold_{NL}$ in a nested fashion. This is in contrast to the solution of problem \eqref{eq:generalPb} by the Simultaneous ANalysis and Design (SAND) 
approach, where the state variables of \eqref{eq:generalPb} are considered in this case as optimization parameters.

{\it REMARK.} While the computational framework proposed in this paper for accelerating the enforcement of the discretized, parametric, linear PDE constraint 
$\Abold(\mubold) \wbold_{L}(\mubold) - \bbold(\mubold) = \mathbf{0}$ is described in the context of the solution of problem \eqref{eq:generalPb} by the NAND approach, it is equally applicable in the context of the solution 
of this problem by the SAND approach. However, its performance may not scale well with the substantially larger number of optimization parameters that is typically encountered in the latter case.
  
\subsection{Gradient-based optimization}
\label{sec:sensitivities}  
  
Because the discretized PDE constraints introduced above may lead to large-scale systems of equations, a gradient-based method is preferred for the solution of problem \eqref{eq:generalPb}.
To this end, this section reviews the computation of the relevant gradients, in order to keep this paper as self-contained as possible.
  
Let
\begin{equation*}
\label{eq:w_all}
\wbold = \begin{bmatrix} \wbold_{L} \\
\wbold_{NL} 
\end{bmatrix} 
\in \Rbb^{N_{w_{L}}+ N_{w_{NL}}}
\end{equation*}
and
\begin{equation}
\label{eq:PDE_constraint}
	\Rbold(\wbold(\mubold) , \mubold) = \begin{bmatrix} \Rbold_L(\wbold_L(\mubold),\mubold)\\ \Rbold_{NL}(\wbold_{NL}(\mubold),\mubold) \end{bmatrix} =\mathbf{0}
\end{equation}
A gradient-based optimization algorithm incurs the computation of the first derivatives of the objective function $f(\wbold(\mubold),\mubold)$, and those of the constraints $c_i(\wbold(\mubold),\mubold)$,
$i=1,~\cdots,~N_c$. Let $q$ denote the generic representation of $f$, any constraint $c_i$, and any QoI whose derivatives with respect to the components of the parameter point $\mubold \in \mathcal D \subset {\mathbb R}^{N_{\mathcal D}}$ 
must be computed. Using the chain rule, the following relations can be derived
\begin{equation}\label{eq:Funcsens}
\frac{dq}{d\mu_i}\left(\wbold(\mubold),\mubold\right) = \frac{\partial q}{\partial \mu_i} \left(\wbold(\mubold),\mubold\right) + \frac{\partial q}{\partial \wbold}\left(\wbold(\mubold),\mubold\right)
\frac{\partial\wbold}{\partial\mu_i}(\mubold), \quad i = 1,~\cdots,~N_{\mathcal D}
\end{equation}
For a given $q$, the computation of the partial derivatives 
$\displaystyle{\frac{\partial q}{\partial \wbold}\left(\wbold(\mubold),\mubold\right)}$ and $\displaystyle{\frac{\partial q}{\partial \mu_i} \left(\wbold(\mubold),\mubold\right)}$ is straightforward.
On the other hand, the computation of $\displaystyle{\frac{\partial\wbold}{\partial\mu_i}(\mubold)}$ warrants some attention. 

The differentiation of the discretized, parametric, PDE constraint \eqref{eq:PDE_constraint} with respect to any parameter component
$\mu_i$ can be written as
\begin{equation*}\label{eq:LinPDEsens}
\frac{d \Rbold}{d \mu_i} \left(\wbold(\mubold),\mubold\right) = \frac{\partial\Rbold}{\partial\mu_i}\left(\wbold(\mubold),\mubold\right) 
+ \frac{\partial\Rbold}{\partial\wbold} \frac{ \partial \wbold }{ \partial \mu_i }\left(\wbold(\mubold),\mubold\right) = \mathbf{0}, \quad i=1,~\cdots,~N_{\mathcal D}
\end{equation*}
which leads to
\begin{equation}\label{eq:REL_LinPDEsens}
\frac{\partial \wbold}{\partial \mu_i }\left(\wbold(\mubold),\mubold\right) = - \left[ \frac{\partial\Rbold}{\partial\wbold} \right]^{-1} \frac{\partial\Rbold}{\partial \mu_i } 
\left(\wbold(\mubold),\mubold\right), \quad i=1,~\cdots,~N_{\mathcal D}
\end{equation}

The differentiation of the discretized, parametric, linear constraint (\ref{eq:LinearForm}) with respect to any parameter 
component $\mu_i$ leads to the following linear system of equations governing
the sensitivities of $\wbold$
\begin{equation*}\label{eq:state_sens}
	\Abold(\mubold) \frac{\partial \wbold}{\partial \mu_i} (\mubold) = \frac{\partial\bbold}{\partial\mu_i}(\mubold) - \frac{\partial\Abold}{\partial\mu_i}(\mubold)\wbold(\mubold), \quad i = 1,~\cdots,~N_{\mathcal D}
\end{equation*}
These state vector sensitivities can be used to compute the gradients of the objective function and constraints. Indeed, substituting the solution of \eqref{eq:REL_LinPDEsens} into \eqref{eq:Funcsens} 
leads to 
\begin{equation}\label{eq:dwdmu_linsyst1}
\frac{dq}{d\mu_i} = \frac{\partial q}{\partial \mu_i}
- \frac{\partial q}{\partial \wbold} \left( \left[ \frac{\partial\Rbold}{\partial\wbold} \right]^{-1} \frac{\partial\Rbold}{\partial \mu_i } \right), \quad i = 1,~\cdots,~N_{\mathcal D}
\end{equation}
which can also be written as
\begin{equation}\label{eq:dwdmu_linsyst2}
\frac{dq}{d\mu_i} = \frac{\partial q}{\partial \mu_i}
- \left[  \frac{\partial\Rbold}{\partial\wbold}^{-T} \frac{\partial q}{\partial \wbold}^T  \right]^T \frac{\partial\Rbold}{\partial \mu_i }, \quad i = 1,~\cdots,~N_{\mathcal D}
\end{equation}
where the superscript denotes here and in the remainder of this paper the transpose operation.
Equations \eqref{eq:dwdmu_linsyst1} and \eqref{eq:dwdmu_linsyst2} describe the {\it direct} and {\it adjoint} approaches for computing the set of sensitivities 
$\displaystyle{\left\{\frac{dq}{d\mu_i}(\wbold(\mubold),\mubold)\right\}_{i=1}^{N_{\mathcal D}}}$, respectively:
\begin{enumerate}
\item In the direct approach, the state sensitivities $\displaystyle{\frac{\partial\wbold}{\partial \mu_i}(\mubold)}$ are first computed by solving the linear system of equations \eqref{eq:REL_LinPDEsens} for each parameter component $\mu_i$, 
	then the sensitivities $\displaystyle{\frac{dq}{d\mu_i}}$ are evaluated using (\ref{eq:Funcsens}).
\item In the adjoint approach, the adjoint vector $\lambdabold_q(\mubold)$ is first computed by solving the linear system of equations 
	$ \displaystyle{\frac{\partial \Rbold}{\partial \wbold}^T \lambdabold_q(\mubold) = \frac{\partial q}{\partial \wbold}\left(\wbold(\mubold),\mubold\right)}$ for $\lambdabold_q(\mubold)$, then all 
	sensitivities $\displaystyle{\frac{dq}{d\mu_i}}$ are computed using \eqref{eq:dwdmu_linsyst2}.
\end{enumerate}
The direct and adjoint approaches require $N_{\mathcal D}$ and $N_c+1$ ``solves'', respectively.  Hence, if $N_{\mathcal D}\leq1+N_c$, the direct approach is preferrable; otherwise, the adjoint approach is preferrable.  

Using the gradients computed above, problem \eqref{eq:generalPb} can be solved, for example, by a nonlinear optimization algorithm such as Sequential Quadratic Programming (SQP) \cite{gill2002} equipped with a 
quasi-Newton approximation of the Hessian matrix, the trust-region method \cite{conn2000trust}, or the interior-point method \cite{wachter2006implementation}. In either case, the treatment of the discretized
PDE constraints is computationally intensive. For example, computing the sensitivities of the linear PDE's state vector with respect to the optimization parameters requires the solution of the linear system of 
equations \eqref{eq:LinearForm} and that of another, large-scale linear system of equations of size $N_{w_L}$. To significantly reduce this cost, model reduction can be applied to \eqref{eq:LinearForm}, and 
reduced-order versions of $\displaystyle{\frac{\partial\bbold}{\partial\mu_i}(\mubold)}$ and $\displaystyle{\frac{\partial\Abold}{\partial\mu_i}(\mubold)}$, $i = 1,~\cdots,~N_{\mathcal D}$, can be computed instead of 
their high-dimensional counterparts.
    
\section{Parametric model order reduction }\label{sec:PROM}

Again, the focus of this paper is on accelerating the solution of the optimization problem \eqref{eq:generalPb}, in the case where the discretized, parametric, PDE constraint \eqref{eq:LinearForm} is 
high-dimensional. In this case, solving the optimization problem \eqref{eq:generalPb} can be prohibitively expensive. For this purpose, model reduction is applied to the discrete system of linear constraints 
in \eqref{eq:PDE_constraint} in order to reduce its dimension, which reduces the cost of enforcing these constraints and therefore reduces the cost of solving problem \eqref{eq:generalPb}.

\subsection{Projection-based model order reduction and representation}

Projection-based model order reduction reduces the dimension of a $\mubold$-parametric HDM such as (\ref{eq:LinearForm}), which governs the state vector $\wbold_L$, by performing the subspace approximation
  \begin{equation}\label{eq:subspaceAssumption}
	  \wbold_L(\mubold) \approx \Vbold(\mubold) \wbold_{L_r}(\mubold)
  \end{equation}
where $\Vbold(\mubold) \in \Rbb^{N_{w_L} \times N_{w_{r}}}$ is a {\it right} ROB of dimension $N_{w_{r}} \ll N_{w_L}$, and $\wbold_{L_r} \in \Rbb^{N_{w_{L_r}}}$ is a vector of generalized coordinates referred to
as the reduced-order state vector. The notation $\Vbold(\mubold)$ emphasizes that $\Vbold$ is a pointwise ROB -- that is, a ROB constructed at the parameter point $\mubold \in \mathcal D$ -- rather than a global ROB constructed to deliver accurate subspace approximations at any parameter point in the design parameter space. 
Typically, $\Vbold(\mubold)$ is constructed using the Proper Orthogonal 
Decomposition (POD) method of snapshots \cite{sirovich1987turbulence}, which can be summarized as the collection of solution snapshots into a matrix, followed by the compression of this matrix using 
the Singular Value Decomposition (SVD) algorithm (see Algorithm \ref{alg:POD}).

\begin{algorithm}
\caption{Proper Orthogonal Decomposition.}
\label{alg:POD}
\begin{algorithmic}[1]
\INPUT Snapshot matrix $\Xbold(\mubold) \in R^{N_{w} \times N_s}$; desired ROB dimension $N_{w_{r}}$
\OUTPUT ROB $\Vbold(\mubold)$
\STATE Compute the thin SVD of $\Xbold(\mubold)$: $\Xbold(\mubold) = \Ubold \Sigmabold \Dbold$, where $\Ubold = \left[ \ubold_1 \,\,\,\, \ubold_2 \,\,\,\, \cdots \,\,\,\, \ubold_{N_s}  \right]$
\STATE $\Vbold(\mubold) = \left[ \ubold_1 \,\,\,\, \ubold_2  \,\,\,\, \cdots \,\,\,\, \ubold_{N_{w_{r}}}  \right]$
\end{algorithmic}
\end{algorithm}
  
Substituting \eqref{eq:subspaceAssumption} into \eqref{eq:LinearForm} and pre-multiplying the resulting system of equations by $\Wbold(\mubold)^T$, where $\Wbold(\mubold) \in \Rbb^{N_{w_L}\times N_{w_{L_r}}}$ is
known as the {\it left} ROB associated with the right ROB $\Vbold(\mubold)$, transforms this system into a square counterpart and leads to the parametric, linear, Petrov-Galerkin PROM
\begin{equation}\label{eq:PROM}
\Wbold(\mubold)^T \Abold(\mubold)\Vbold(\mubold) \wbold_{L_r} = \Wbold(\mubold)^T \bbold(\mubold)
\end{equation}
This PROM can be described by the low-dimensional matrix $\Abold_r(\mubold) = \Wbold(\mubold)^T \Abold(\mubold) \Vbold(\mubold) \in \Rbb^{N_{w_{L_r}} \times N_{w_{L_r}}}$, the low-dimensional vector 
$\bbold_r(\mubold) = \Wbold(\mubold)^T \bbold(\mubold) \in \Rbb^{N_{w_{L_r}}}$, and therefore the dublet $\left(\Abold_r(\mubold), \bbold_r(\mubold)\right)$. 
If $\Wbold(\mubold) = \Vbold(\mubold)$, this PROM becomes a Galerkin PROM.

The sensitivities of the reduced-order state vector $\wbold_{L_r}$ with respect to the parameter components 
$\mu_i$, $i = 1,~\cdots,~N_{\mathcal D}$, can be computed by following the same approach described in 
Section \ref{sec:sensitivities}. In this case, this approach entails the computation of the reduced-order sensitivity matrices and sensitivity vectors 
$\displaystyle{\frac{\partial \Abold_r}{\partial \mu_i}(\mubold) \in \Rbb^{N_{w_{L_r}} \times N_{w_{L_r}}}} $ and $\displaystyle{\frac{\partial \bbold_r}{\partial \mu_i}(\mubold) \in \Rbb^{N_{w_{L_r}}} }$,
respectively, $i = 1,~\cdots,~N_{\mathcal D}$. It follows that the PROM (\ref{eq:PROM}) and its associated reduced-order sensitivity equations can be collectively represented by the tuple of 
low-dimensional operators 
\begin{equation}\label{eq:tuple}
	\mathcal{A}_r(\mubold) = \left\{ \Abold_r(\mubold), \bbold_r(\mubold), \frac{\partial \Abold_r}{\partial \mu_i}(\mubold), \frac{\partial \bbold_r}{\partial \mu_i}(\mubold), \quad
	i = 1,~\cdots,~N_{\mathcal D} \right\}
\end{equation}  

In summary, the subspace approximation (\ref{eq:subspaceAssumption}) transforms the original constrained optimization problem (\ref{eq:generalPb}) into the PROM-constrained optimization problem
\begin{equation}\label{eq:ROMPb}
\boxed{
\begin{aligned}
& \underset{\wbold\in\mathcal{W},\mubold\in\mathcal{D}}{\text{min}}
& & f(\Vbold(\mubold) \wbold_{L_r}(\mubold),\wbold_{NL}(\mubold),\mubold) \\
& \text{s.t.} & &  \cbold\left(\Vbold(\mubold)\wbold_{L_r}(\mubold),\wbold_{NL}(\mubold),\mubold)\right) \leq \boldsymbol{0} \\
& & &  \Abold_r(\mubold) \wbold_{L_r}(\mubold) = \bbold_r(\mubold) \\
& & &  \Rbold_{NL}(\wbold_{NL}(\mubold),\mubold)= \mathbf{0} \\
 \end{aligned}
}
\end{equation}
where the discretized, parametric, linear PDE constraint and its associated sensitivity equations are low-dimensional. Just like the original constrained optimization problem (\ref{eq:generalPb}),
problem (\ref{eq:ROMPb}) can be solved using any preferred, gradient-based optimization algorithm.
      
\subsection{Concept of a database of pointwise linear PROMs and associated sensitivities}

As stated in Section \ref{sec:Introduction}, a PROM constructed at a parameter point $\mubold \in \mathcal D$ does not necessarily perform well at a different parameter point $\mubold^{\prime}$ in this space. Addressing this issue by constructing a new PROM at each newly queried parameter point $\mubold^{\prime}$ is computationally inefficient, as it entails: the computation of new solution snapshots 
at the parameter point $\mubold^{\prime}$; the compression of these snapshots using, for example, the SVD algorithm in order to compute a new right ROB $\Vbold(\mubold^{\prime})$; and performing additional 
matrix-matrix and matrix-vector computations to compute, in the context of this paper, the tuple $\mathcal{A}_r(\mubold^{\prime}) = \displaystyle{\left\{ \Abold_r(\mubold^{\prime}), \bbold_r(\mubold^{\prime}), \frac{\partial \Abold_r}{\partial \mu_i}(\mubold^{\prime}), \frac{\partial \bbold_r}{\partial \mu_i}(\mubold^{\prime}), \quad i = 1,~\cdots,~N_{\mathcal D}\right\}}$ representing the linear PROM at the parameter point 
$\mubold^{\prime}$.

An alternative methodology for addressing efficiently the parameter dependence of a linear PROM based on the concept of a database of 
reduced-order information was proposed and developed in 
\cite{amsallem2010towards}, and more recently refined in \cite{amsallem2016realtime}. It consists in: constructing offline a database of pointwise, linear PROMs; 
equipping this database with a family of algorithms for interpolation on matrix manifolds \cite{amsallem2008interpolation, amsallem2009method, amsallem2011online}; and applying these algorithms online 
to build in real time a PROM at a queried but unsampled parameter point 
$\tilde \mubold \in \mathcal D \subset {\mathbb R}^{N_{\mathcal D}}$. 
In the context of the parametric, linear PROM (\ref{eq:PROM}), this offline-online methodology can be described as follows:
\begin{enumerate} 
	\item {\bf Offline} 
		\begin{itemize} 
			\item Sample $N_{\mathcal{DB}}$ parameter points $\mubold^j$ in the given design parameter space $\mathcal D$.  
			\item At each sampled parameter point $\mubold^j, j = 1, \cdots, N_{\mathcal{DB}}$, construct a pointwise, linear PROM 
				$\mathcal{A}_r(\mubold^j) = \displaystyle{\left\{ \Abold_r(\mubold^j), \bbold_r(\mubold^j), \frac{\partial \Abold_r}{\partial \mu_i}(\mubold^j), \frac{\partial \bbold_r}{\partial \mu_i}(\mubold^j), \quad i = 1,~\cdots,~N_{\mathcal D}\right\}}$.
			\item Store the set of pre-computed linear PROMs in a database $\mathcal{DB}$ (in practice, the derivatives $\displaystyle{\frac{\partial \Abold_r}{\partial \mu_i}(\mubold^j)}$ 
			and $\displaystyle{\frac{\partial \bbold_r}{\partial \mu_i}(\mubold^j)}$ need not be stored, because they can be computed online as explained in Section~\ref{sec:interpolation}).
		\end{itemize}  
	\item {\bf Online} 

		For each queried but unsampled parameter point $\tilde \mubold$, construct a pointwise, linear PROM $\mathcal{A}_r(\tilde \mubold) = \displaystyle{\left\{ \Abold_r(\tilde \mubold), \bbold_r(\tilde \mubold), \frac{\partial \Abold_r}{\partial \mu_i}(\tilde\mubold), \frac{\partial \bbold_r}{\partial \mu_i}(\tilde\mubold), \quad i = 1,~\cdots,~N_{\mathcal D}\right\}}$ as follows:
			\begin{itemize} 
				\item For each reduced-order quantity in $\mathcal{A}_r(\mubold)$, identify {\it \`a priori} a matrix manifold where it lies, for example, by identifying the most important 
					algebraic properties characterizing this reduced-order quantity.
				\item Interpolate in real time the pre-computed reduced-order quantities on their identified manifolds, at the queried but unsampled parameter point $\tilde\mubold$.
			\end{itemize} 
	\end{enumerate}

In order to extend the scope of this concept of a database of pointwise, linear PROMs to the real-time solution of the constraint equation (\ref{eq:LinearForm}), in view of accelerating the solution of 
the optimization problem \eqref{eq:generalPb},  the state of the art of this concept is next advanced as follows. First, a computationally efficient parameter sampling procedure based on a saturation assumption
is described. Then, a novel algorithm for interpolating in real time the sensitivities of a pointwise, linear PROM with respect to parameter variations is presented. 

\section{Parameter sampling using a residual-based error indicator}
\label{sec:database}

The typical construction of an $N_{\mathcal{DB}}$-point database 
$\mathcal{DB} = \left\{ \mubold^j, \mathcal{A}_r(\mubold^j)\right\}_{j=1}^{N_\mathcal{DB}}$ for the 
purpose of supporting the concept of interpolation on matrix manifolds, which addresses the parameter dependence of a linear PROM, 
is guided by a parameter sampling procedure that selects the parameter points in the design parameter space $\mathcal D$
where to pre-compute linear PROMs of interest. Specifically, such a procedure should lead to a database $\mathcal{DB}$ that has the 
smallest possible size $N_{\mathcal {DB}}$, in order to minmize as much as possible the computational cost of the offline stage. 
It should also populate this database in a manner that delivers interpolation accuracy in the entirety of $\mathcal D$ -- that is,
$\forall \mubold \in \mathcal D$. 

Parameter sampling procedures can be classified in two categories: non-adaptive procedures that sample the given design parameter space
{\it a priori}; and adaptive, error-informed procedures that incorporate a feedback in the sampling process. Non-adaptive parameter 
sampling procedures are typically easier to implement. Adaptive counterparts however can be expected to perform better at producing 
databases $\mathcal{DB}$ that satisfy the desiderata outlined above.

The simplest non-adaptive parameter sampling procedure is perhaps the so-called Full Factorial (FF) sampling procedure, which samples 
the given design parameter space on an $N_{\mathcal D}$-dimensional grid. In principle, this grid can be made sufficiently fine
to deliver accurate approximations such as those
to be performed here using interpolation on a matrix manifold. However, the FF sampling procedure rapidly becomes 
unfeasible, as the number of parameter points it requires for achieving accuracy grows exponentially with $N_{\mathcal {DB}}$.
The alternative Latin Hypercube Sampling (LHS) procedure, which is a non-adaptive but sparse procedure, scales only linearly with 
$N_{\mathcal {DB}}$. However, because it has no explicit awareness of where the approximated solution will be inaccurate, it typically 
leads to oversampling when attempting to meet the aforementioned accuracy desideratum. 

It follows that for high-dimensional design parameter spaces, an adaptive approach where the parameter points are sampled in an incremental
manner and the location of the next sampled parameter point is identified by drawing information from the already sampled parameter points 
is often necessary. Such an approach is typically iterative and referred to as a greedy procedure. It initializes the construction 
of the database 
by populating it with one or a small number of parameter points located, for example, in the center of $\mathcal D$ and/or on its 
boundaries. At each iteration, it considers a set of candidate parameter points for sampling and selects that parameter point where some error estimator
or indicator $e(\mubold, \mathcal{DB})$ is maximized \cite{amsallem2013design, 
bui2008model, bui2008parametric, grepl2005posteriori, hesthaven2014efficient, paul2014adaptive,veroy2003posteriori}.  
After $n$ greedy iterations, it produces a database $\mathcal{DB}_n$ that can be characterized as follows
\begin{equation*}\label{eq:greedy}
\mathcal{DB}_n = \left\{  \mathcal{DB}_{n-1}, \mubold^n, \mathcal{A}_r(\mubold^n) \right\} \quad \text{ where } \quad \mubold^n =  \arg \max_{\mubold \in \Xi} e(\mubold, \mathcal{DB}_{n-1})
\end{equation*}
where $\Xi = \{ \mubold_c^1, \cdots, \mubold_c^{N_\Xi} \} $ denotes the set of $N_\Xi$ candidate parameter points in the design 
parameter space $\mathcal D$. This set can be fixed throughout the iterations or refreshed at each iteration, but in any case, 
should be chosen large enough to faithfully represent $\mathcal{D}$. The iterations performed by such a greedy 
procedure are terminated when the error estimator or indicator falls below a specified threshold, or the pre-set maximum size of 
the database or maximum run-time allocated for the execution of the procedure has been reached. 

The most pertinent error for an adaptive sampling procedure is the true error 
$e(\mubold_c^j) = \Vert\wbold_{L}(\mubold_c^j) - \wbold_L^I(\mubold_c^j)\Vert| = 
|\Vert\wbold_{L}(\mubold_c^j) - \widetilde \Vbold(\mubold_c^j)\wbold_{L_r}^I(\mubold_c^j)\Vert$, 
where $\wbold_L(\mubold_c^j)$ and $\wbold_{L_r}^I(\mubold_c^j)$ denote the solutions computed at the candidate parameter point 
$\mubold_c^j \in \mathcal D$ of the HDM and PROM problems constructed and interpolated at $\mubold_c^j$, respectively, and
$\widetilde \Vbold (\mubold_c^j)$ is an approximation of $\Vbold (\mubold_c^j)$ that can be computed in real time (see Section
\ref{sec:RBEI}). However, whereas interpolating the PROMs, solving the reduced-order problems, and reconstructing the corresponding 
high-dimensional solutions at the set of candidate parameter points is feasible, because each of these operations can be performed
relatively fast if not in real time, the direct computation of the high-dimensional solutions for the entire set $\Xi$ is typically 
prohibitively expensive. Consequently, the true error is replaced in practice by: an estimation of this error that is based on a 
proven error bound, when such a bound is available; or by an error indicator otherwise. Error bounds are known for the solutions
of some relatively simple elliptic \cite{prud2002reliable, veroy2003posteriori}, parabolic \cite{grepl2005posteriori}, and 
hyperbolic \cite{haasdonk2008reduced} PDEs, as well as for the solutions of some simple linear time-invariant systems 
\cite{haasdonk2011efficient}. However, these bounds typically require the computation of an inf-sup constant, which
is seldom a simple task, and are not usually tight \cite{amsallem2014posteriori}. For this reason, computationally feasible error indicators based on some
norm of a relevant residual evaluated using the reconstructed high-dimensional solution are used instead \cite{bui2008model, bui2008parametric}. In the context of this paper, 
such an error indicator can be expressed as $\Vert\Rbold_L\left(\widetilde\Vbold(\mubold_c^j)\wbold_{L_r}^I(\mubold_c^j)\right)\Vert$, where $\Rbold_L$ is defined in \eqref{eq:generalPb}. 
Nevertheless, even the computation of such error indicators can rapidly become overwhelming when
$N_\Xi$ becomes very large, for example, because the dimension $N_{\mathcal D}$ of the design parameter space $\mathcal D$ is increased.
One approach for addressing this issue is to reduce as much as possible the number of parameter points that must be sampled
to achieve the desired accuracy -- that is, to minimize as much as possible the size $N_{\mathcal {DB}}$ of the database 
$\mathcal {DB}$. To this end, an accelerated parameter sampling procedure based on an appropriate saturation assumption is proposed next.

\subsection{Accelerated parameter sampling procedure based on a saturation assumption}
\label{sec:EGP}

Here, the set of $N_\Xi$ candidate parameter points in the design parameter space $\mathcal D$, $\Xi = \{ \mubold_c^1, \cdots, \mubold_c^{N_\Xi} \}$, is fixed throughout all iterations
of the greedy procedure. In this case, this parameter sampling procedure can be further accelerated by using a saturation technique where previous evaluations of the error indicator are used to 
prevent, whenever possible, unnecessary additional evaluations of this indicator. This is possible if the accuracy of the database can be assumed to increase only to a certain extent as more parameter 
points are sampled, as this assumption makes it possible to omit the evaluation of the error indicator at some parameter points.

The incorporation of a saturation assumption in a greedy procedure was first proposed in \cite{hesthaven2014efficient}. Here, the efficiency of this technique is further improved by applying it at
each $m$-th iteration of the greedy procedure not to $\Xi$, but to a random subset $\Pi_{m}$ of $\Xi$. Specifically, a fixed number of $N_{\Pi}$ parameter points are 
randomly selected in $\Xi$ at each greedy iteration to assess the accuracy of the database being constructed, and a saturation constant is defined as follows.

\begin{mydef} \textbf{Saturation Constant}\\
	Let $e(\mubold; \mathcal{DB}_{m})$ denote an error indicator depending on the design parameter vector $\mubold$ and the database $\mathcal{DB}_{m}$. Let also $\{\mathcal{DB}_{m}\}_{m=1}^n$
	denote a sequence of nested databases -- that is, a set of databases satisfying $\mathcal{DB}_{m} \subset \mathcal{DB}_{n}$ for all $1\leq m<n$. The saturation constant $\tau_s > 0$ is defined 
	as 
\begin{equation}\label{eq:SaturationConstantDefinition}
\tau_s = \underset{1\leq m<n, \, \mubold\in\mathcal{D}}{\text{sup }}\hspace{7pt} \frac{e(\mubold;\mathcal{DB}_{n})}{e(\mubold; \mathcal{DB}_{m})}
\end{equation} 
\end{mydef}

From the above definition, it follows that: $e(\mubold;\mathcal{DB}_{n}) \leq \tau_s e(\mubold; \mathcal{DB}_{m})$, for all $1\leq m<n$;
and $\tau_s < 1$ implies that $e(\mubold;\mathcal{DB}_{n}) < e(\mubold; \mathcal{DB}_{m})$ for all $1\leq m<n$. In other words, $e(\mubold;\cdot)$ strictly decreases as more parameter points are 
sampled during the construction of the database. Similarly, $\tau_s = 1$ implies a monotone decrease in $e(\mubold;\cdot)$ as the number of sampled parameter points is increased during
the construction of the database. If $\tau_s > 1$, $e(\mubold;\cdot)$ may increase at some parameter point $\mubold \in \mathcal{DB}$ for certain iterations of the greedy procedure. 

Now, consider the additional two definitions:

\begin{itemize}
\item Let $e_{\text{profile}}(\mubold_c)$ denote the most recently computed error indicator at the candidate parameter point $\mubold_c\in \Pi_{m} \subset \Xi $, where $\Pi_{m}$ is
	a subset of $\Xi$ randomly selected during the $m$-th iteration of the greedy procedure. Note that because of the randomness aspect of the construction of the subset of candidate parameter
		points $\Pi_{m}$, $\mubold_c$ might not have been selected at iteration $m-1$ for error indication. Hence in general, $e_{\text{profile}}(\mubold_c) \ne e(\mubold_c; \mathcal{DB}_{m-1})$.
		Consequently, if $\mubold_c$ was not previously selected for computing the error indicator, $e_{\text{profile}}(\mubold_c)$ is set to $\infty$. 
\item The following definition is associated with a subiterative process within the $m$-th iteration of the greedy procedure. Let $e^{(l)}_{\max}(\mathcal{DB}_{m})$ denote the current maximum indicated error at subiteration $l$ 
of iteration $m$ of the greedy procedure, where one error indication is computed at a candidate parameter point in $\Pi_m \subset \Xi$. This iterative maximum indicated error is computed as follows
\begin{equation*}  
e^{(l)}_{\max}(\mathcal{DB}_{m}) = 
\begin{cases}
0                                                                                        & \text{for } l = 1 (\mbox{initialization}) \\
\underset{j= k_1,\ldots,k_{l-1}}{\text{max }}\hspace{3pt} e(\mubold_c^j; \mathcal{DB}_{m}) & \text{for } 1 < l \leq N_{\Pi} \text{ and } \mubold_c^j \in \Pi_{m}
\end{cases}
\end{equation*} 
\end{itemize}  
where $l$ also denotes a local indexing of a candidate parameter point within $\Pi_m$, and $k_l$ denotes the corresponding global indexing of this candidate parameter point within $\Xi \supset \Pi_m$.
           
Based on the two above definitions and that of the saturation constant~(\ref{eq:SaturationConstantDefinition}), it follows that if $\tau_s e_{\text{profile}}(\mubold_c^i) < e^{(l)}_{\max}(\mathcal{DB}_{m})$, then 
$e(\mubold_c^i; \mathcal{DB}_{m})$ is guaranteed to be less than $e^{(l)}_{\max}(\mathcal{DB}_{m})$ prior to the evaluation of the error indicator at $\mubold_c^i$. This is because 
$e(\mubold_c^i; \mathcal{DB}_{m}) \leq \tau_s e_{\text{profile}}(\mubold_c^i)$, due to the definition of the saturation constant. Hence, it is not necessary to evaluate $e(\mubold_c^i; \mathcal{DB}_{m})$
if $\tau_s e_{\text{profile}}(\mubold_c^i) < e^{(l)}_{\max}(\mathcal{DB}_{m})$. This implies that the simple inequality check between $\tau_s e_{\text{profile}}(\mubold_c^i)$ and 
$e^{(l)}_{\max}(\mathcal{DB}_{m})$, whose values are readily available, can avoid some unnecessary evaluations of the error indicator. Then, convergence of the greedy procedure can be monitored
as follows. A check is performed to verify whether $e_{\max}^{(l)}(\mathcal{DB}_{m}) < \epsilon_{tol}$ or not, where $\epsilon_{tol}$ is a specified threshold tolerance for assessing the accuracy 
of the constructed database. For this purpose, the error indicator is evaluated at all parameter points of a larger subset of $\Xi$ of size $N_{\Pi}^{\prime} > N_{\Pi}$ that is also randomly selected. 
If $e_{\max}^{(l)}(\mathcal{DB}_{m}) < \epsilon_{tol}$ holds for this larger subset of candidate parameters, the greedy procedure is terminated. 

In general, the value of $\tau_s$ is set to 1 or 2: 1 for an aggressive approach; and 2 for a conservative approach. 

Algorithm \ref{al:ROMinterpolationgreedyalgorithm} summarizes the greedy procedure for parameter sampling equipped with the saturation technique.  

\label{sec:saturationgreedy} 
\begin{algorithm}[!h]
\caption{Greedy procedure incorporating a randomized saturation assumption.}
\label{al:ROMinterpolationgreedyalgorithm}
\textbf{Input:} A set of candidate parameter points $\Xi\subset\mathcal{D}$, its size $N_\Xi$, a tolerance $\epsilon_{tol} > 0$, a saturation constant $\tau_s$, and two sizes $N_\Pi$
	and $N_{\Pi}^{\prime} > N_{\Pi}$ for the two subsets of $\Xi$ to be randomly selected at each greedy iteration\\
       \textbf{Output:} PROM database $\mathcal{DB}_{m}$
         \begin{algorithmic}[1]
           \STATE Choose an initial parameter point $\mubold_c^1\in\Xi$ and compute $ \mathcal{A}_r(\mubold_c^1)$
           \STATE Set $\mathcal{DB}_1 = \{ (\mubold_c^1, \mathcal{A}_r(\mubold_c^1)) \}$
           \STATE Set $e_{\text{profile}}(\mubold) = \infty$ for all $\mubold \in \Xi$
           \WHILE { TRUE }
              \STATE Generate a random subset $\Pi_{m}$ of $\Xi$ of dimension $N_\Pi$
              \STATE Set $e_{\max}^{(l)} = 0$ 
              \FOR { $\mubold_c^{k_l} \in \Pi_{m}, l=1,\ldots,N_\Pi$ } 
                 \IF {$\tau_s e_{\text{profile}}(\mubold_c^{k_l}) > e_{\max}^{(l)}$}
                    \STATE Compute $e(\mubold_c^{k_l};\mathcal{DB}_{m})$ and set 
                $e_{\text{profile}}(\mubold_c^{k_l}) = e(\mubold_c^{k_l};\mathcal{DB}_{m})$
                    \IF { $e(\mubold_c^{k_l};\mathcal{DB}_{m}) > e_{\max}^{(l)}$ }
                      \STATE Set $e_{\max}^{(l)}  = e(\mubold_c^{k_l};\mathcal{DB}_{m})$ and $\mubold_{m+1} = \mubold_c^{k_l}$
                    \ENDIF 
                 \ENDIF 
              \ENDFOR
              \IF { $e_{\max}^{(l)} < \epsilon_{tol}$ }
		 \STATE Check convergence using the larger subset of candidate parameter points of size $N_{\Pi}^{\prime} > N_{\Pi}$
                 \IF {convergence is reached} 
                    \STATE Terminate the parameter sampling procedure
                 \ENDIF
              \ENDIF
              \STATE Compute $\mathcal{A}_r(\mubold_c^{m+1})$
		 \STATE Set $\mathcal{DB}_{m+1} = \mathcal{DB}_{m} \cup \{ (\mubold_c^{m+1}, \mathcal{A}_r(\mubold_c^{m+1}) )\}$ 
              \STATE $m \gets m+1$
           \ENDWHILE
         \end{algorithmic}
       \end{algorithm}

\subsection{Computation of the residual-based error indicator}
\label{sec:RBEI}

To complete the description of the greedy parameter sampling procedure proposed in this paper for guiding the construction of a database of linear, pointwise PROMs, it remains to discuss one subtlety 
regarding the practical computation of a residual-based error indicator such as
\begin{equation}\label{eq:residual}
\left\Vert\Rbold_L\left(\Vbold(\mubold_c^j)\wbold_{L_r}^I(\mubold_c^j)\right)\right\Vert
\end{equation}
where all quantities appearing in \eqref{eq:residual} have been defined in the introduction part of Section \ref{sec:database}. 
Note that neither the right ROB $\Vbold(\mubold)$ nor the left counterpart $\Wbold(\mubold)$ is part of the definition of the tuple 
\break $\mathcal{A}_r(\mubold) = \displaystyle{\left\{ \Abold_r(\mubold), \bbold_r(\mubold), \frac{\partial \Abold_r}{\partial \mu_i}(\mubold), \frac{\partial \bbold_r}{\partial \mu_i}(\mubold), \quad i = 1,~\cdots,~N_{\mathcal D}\right\}}$ representing the linear 
PROM constructed or reconstructed at the parameter point $\mubold$.
Note also that the computational complexity of the construction of a ROB such as $\Vbold(\mubold)$, for $\mubold \in \mathcal D$,
scales with the dimension of the HDM:
therefore, the construction of such a ROB at each candidate parameter point $\mubold_c^j \in \Xi$ is not feasible. 
Instead, a reasonable
approximation $\widetilde \Vbold(\mubold_c^j)$ of $\Vbold(\mubold_c^j)$ is sought-after here. For this purpose, computing 
$\widetilde \Vbold(\mubold_c^j)$ by interpolating on the Grassman manifold \cite{amsallem2008interpolation} the pre-computed ROBs at 
the previously sampled parameter points, which could be stored in this case in the database $\mathcal {DB}$, is feasible, has
been shown in \cite{amsallem2008interpolation} to deliver a good level of accuracy even in the context of a sparsly populated
database $\mathcal DB$, but unfortunately requires the solution of an auxiliary problem whose computational complexity scales with 
the dimension of the HDM (see \cite{amsallem2008interpolation} for details). Hence, given that in the context of a greedy 
procedure the level of accuracy of an approximation is not paramount for the selection among the elements of $\Xi$ of the best
parameter point to be sampled, it is proposed here to store the aforementioned pre-computed ROBs in $\mathcal {DB}$ and reconstruct 
$\widetilde \Vbold(\mubold_c^j)$ for each $\mubold_c^j \in \Xi$ using a constant extrapolation procedure based on the nearest
sampled parameter point $\mubold^j$. This approach for computing $\widetilde \Vbold(\mubold_c^j)$, for each $\mubold_c^j \in \Xi$,
leads to a procedure for constructing the database of interest $\mathcal {DB}$ that can be summarized as follows
\begin{eqnarray*}\label{eq:basisDB}
	\mathcal{DB} &=& \left\{ \mubold^j , \mathcal{A}_r(\mubold^j), \Vbold(\mubold^j) \right\}_{j=1}^{N_{\mathcal{DB}}}\\
	\hbox{and} \, \Rbold_L\left(\Vbold(\mubold_c^j)\wbold_{L_r}^I(\mubold_c^j)\right) &\approx&\Rbold_L\left(\Vbold(\hat{\mubold}) \wbold_{L_r}^I(\mubold_c^j)\right),\quad\hbox{where}\quad
	\hat{\mubold} = \arg \min_{\mubold^{\prime} \in \mathcal{DB}}\Vert\mubold^{\prime} - \mubold\Vert_2^2
\end{eqnarray*}
This approach is not only feasible, but also computationally efficient.

 \begin{figure}
        \begin{center}
          \includegraphics[scale=0.50]{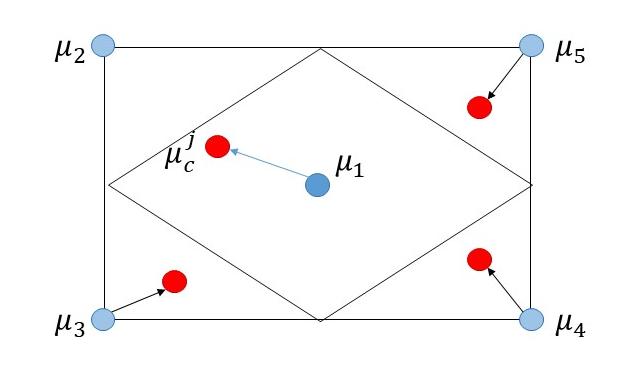}
        \end{center}
	\vspace{-3mm}
        \caption{Greedy parameter sampling procedure: evaluation of the residual-based error indicator using a constant extrapolation
	 approach for approximating $\Vbold_{\mubold_c^j}$.}
        \label{fig:residualBasedEE}
  \end{figure} 

\section{Interpolation on a matrix manifold}
\label{sec:interpolation}

Suppose that an $N_{\mathcal{DB}}$-point database $\mathcal{DB} = \left\{ \mubold^j, \mathcal{A}_r(\mubold^j), \Vbold(\mubold^j) \right\}_{j=1}^{N_\mathcal{DB}}$ is
constructed as described above, for the purpose of accelerating the solution of the reduced-order optimization problem \eqref{eq:ROMPb} by a 
gradient-based algorithm. At some point of the optimization process, such algorithm may query an unsampled parameter point $\tilde {\mubold}\in \mathcal{D}$.
Here, the idea is to construct in real time at $\tilde \mubold$ the tuple $\mathcal{A}_r(\tilde{\mubold})$ representing the linear PROM at this point, 
by interpolating the content of $\mathcal DB$.

To this end, it is first noted that the straightforward (or direct)
interpolation of the operators defining $\mathcal{A}_r(\tilde{\mubold})$ -- that is, the approach
where each entry of each operator is interpolated independently from the others -- is bound to deliver inaccurate results for the following reasons:
\begin{enumerate}
\item Because the PROM operators stored in $\mathcal DB$ may have been constructed using different Petrov-Galerkin bases -- that is, for $i \ne j$,
	$\left(\{\Wbold(\mubold^i), \Vbold(\mubold^i)\}\right) \neq \left(\{\Wbold(\mubold^j), \Vbold(\mubold^j)\}\right)$, the 
		reduced-order operators defining the tuples $\mathcal{A}_r(\mubold^j)$ may have been constructed for different generalized coordinate systems \cite{amsallem2011online}.  
        In this case, a correct interpolation requires first transforming the content of the database $\mathcal DB$ into one that uses a 
        single generalized coordinate system.
\item The direct interpolation approach does not necessarily preserve some important properties of the PROM operators being interpolated, such as
      nonsingularity, orthogonality, or positive-semi-definiteness \cite{amsallem2011online, amsallem2008interpolation}.
\end{enumerate}

Both issues raised above can be addressed by first ``rotating'' the PROM operators defining each tuple $\mathcal{A}_r(\tilde{\mubold})$ in order
to achieve consistency in the generalized coordinate system, then interpolating these operators on appropriate matrix 
manifolds \cite{amsallem2011online}.

\subsection{Interpolation of projection-based reduced-order models}

\subsubsection{Enforcement of consistency}
\label{sec:ROT}

Starting from the observation that any ROB $\Vbold(\mubold)$ remains a basis of the same subspace after multiplication by an orthogonal matrix $\Qbold$,  
the following classes of equivalence can be defined for each subspace $\mathcal{S}(\mubold) = \text{range}(\Vbold(\mubold)\Qbold)$
\begin{equation*}\label{eq:equivalenceClasses}
	c\left(\Vbold(\mubold)\right) = \{ \Vbold(\mubold)\Qbold ~~|~~ \Qbold \in \mathcal{O}(N_{w_{L_r}})  \}
\end{equation*}
where $\mathcal{O}(N_{w_{L_r}})$ denotes the set of orthogonal matrices in $\Rbb^{N_{w_{L_r}} \times N_{w_{L_r}}}$. Each equivalence class defines 
a set of bases leading to PROM operators that act in the same subspace, but use for this purpose different systems of generalized coordinates.

The corresponding classes of equivalence for the PROM operators themselves are 
\begin{equation*}\label{eq:PROMequivalenceClasses}
	c\left(\Abold_r(\mubold), \bbold_r(\mubold), \frac{\partial \Abold_r}{\partial \mu_i}(\mubold), \frac{\partial \bbold_r}{\partial \mu_i}(\mubold)\right) = \left\{ \left. \Qbold^T\Abold_r(\mubold) \Qbold, \Qbold^T \bbold_r(\mubold), \Qbold^T\frac{\partial \Abold_r}{\partial \mu_i}(\mubold)\Qbold, \Qbold^T\frac{\partial \bbold_r}{\partial \mu_i}(\mubold) ~~~ \right| \Qbold \in \mathcal{O}(N_{w_{L_r}})  \right\}
\end{equation*}

Hence, given a reference ROB $\Vbold(\mubold^{\circ})$ and any other ROB of interest $\Vbold(\mubold^j)$, the latter ROB can be transformed into 
an equivalent basis that uses the same system of generalized coordinates as that used by $\Vbold(\mubold^{\circ})$, through its right multiplication by 
the rotation matrix $\Qbold$ that solves the following minimization problem
\begin{equation*}\label{eq:procrustes}
        \begin{aligned}
		& \underset{\Qbold \in \mathcal{O}(N_{w_{L_r}})}{\text{min}}
		& & \left |\left| \Vbold(\mubold^j)\Qbold - \Vbold(\mubold^{\circ}) \right|\right|_F^2
 \end{aligned}
    \end{equation*}
This problem, which is known as the orthogonal Procustes problem, has the following analytical solution
\begin{equation}\label{eq:procrustesSolution}
	\Qbold^{j} = \Ubold^j \Zbold^{j^T}
\end{equation}
where $\Ubold^j$ and $\Zbold^j$ are given by the SVD of the matrix product $\Vbold(\mubold^j)^T\Vbold(\mubold^{\circ})$ -- that is,
\begin{equation*}
	\Vbold(\mubold^j)^T\Vbold(\mubold^{\circ}) =\Ubold^j \Sigmabold^j \Zbold^{j^T}
\end{equation*}

The corresponding transformation of the tuple $\mathcal{A}_r(\mubold^j)$ is given by
\begin{equation}\label{eq:rotatedPROM}
	\mathcal{A}_r^{\star}(\mubold^j) = \left\{\underbrace{\Qbold^{j^T}\Abold_r(\mubold^j) \Qbold^{j}}_{\Abold_r^{\star}(\mubold^j)}, 
	\underbrace{\Qbold^{j^T} \bbold_r(\mubold^j)}_{\bbold_r^{\star}(\mubold^j)}, 
	\underbrace{\Qbold^{j^T}\frac{\partial \Abold_r}{\partial \mu_i}(\mubold^j)\Qbold^j}_{\displaystyle{\frac{\partial \Abold_r^{\star}}{\partial \mu_i}(\mubold^j)}}, 
	\underbrace{\Qbold^{j^T}\frac{\partial \bbold_r}{\partial \mu_i}(\mubold^j)}_{\displaystyle{\frac{\partial \bbold_r^{\star}}{\partial \mu_i}(\mubold^j)}}, \quad i = 1,~\cdots,~N_{\mathcal D}  \right\}
\end{equation}

In summary, given an $N_{\mathcal{DB}}$-point database $\mathcal{DB} = \left\{ \mubold^j, \mathcal{A}_r(\mubold^j), \Vbold(\mubold^j) \right\}_{j=1}^{N_\mathcal{DB}}$ where
all tuples $\mathcal{A}_r(\mubold^j)$ are not necessarily consistent -- in the sense that they are not constructed for the same generalized 
coordinate system -- the following procedure may be applied to transform the content of this database into an equivalent one that is consistent: first, 
identify a reference parameter $\mubold^{\circ}$; then, rotate $ \Vbold(\mubold^j)$ and all PROMs stored in the database using 
\eqref{eq:procrustesSolution} and \eqref{eq:rotatedPROM}.

\subsubsection{Logarithm and exponential maps and interpolation in the tangent space} 
\label{sec:LAEMS}

Let
\begin{equation}\label{eq:rotatedDatabase}
	\mathcal{DB}^{\star} = \left\{ \mubold^j , \mathcal{A}_r^{\star}(\mubold^j), \Vbold^{\star}(\mubold^j) \right\}_{j=1}^{N_{DB}}
\end{equation}
where $\Vbold^{\star}(\mubold^j) = \Vbold(\mubold^j)\Qbold^j$ \big(see (\ref{eq:procrustesSolution})\big), 
be an $N_{\mathcal{DB}}$-point database of consistent, linear PROMs. For each reduced-order operator contributing to the definition 
of a tuple $\mathcal{A}_r^{\star}$ representing a linear PROM -- which in this case, is a reduced-order matrix -- let 
$\mathcal M \subset \Rbb^{N_{w_{L_r}} \times M_{w_{L_r}}}$ denote the most important, known, matrix manifold on which 
this operator lies. For this purpose, any reduced-order vector contributing to the definition of the tuple $\mathcal{A}_r^{\star}$
is viewed here and throughout the remainder of this paper, where appropriate, as a single column matrix (hence, $M_{w_{L_r}} = N_{w_{L_r}}$ or $M_{w_{L_r}} = 1$). 
Such a manifold $\cal M$ can be identified, for example, by considering the most important 
algebraic property characterizing the reduced-order matrix to be interpolated
(for example, symmetric positive definiteness in the case of the 
reduced-order mass matrix). Then, when interpolating this reduced-order matrix at a queried but unsmapled parameter point 
$\tilde \mubold \in \mathcal{D}$ -- which contributes to the interpolation of the linear PROM represented by the tuple 
$\mathcal{A}_r(\tilde \mubold)$ -- this property can be preserved by performing the interpolation on the identified
matrix manifold $\mathcal M$. Such an interpolation is described below for a generic reduced-order matrix of the tuple 
$\mathcal{A}_r^{\star}(\mubold)$, and the manifold of interest $\cal M$ on which it lies.

\begin{figure}[t]
\begin{center}
\includegraphics[scale=0.40]{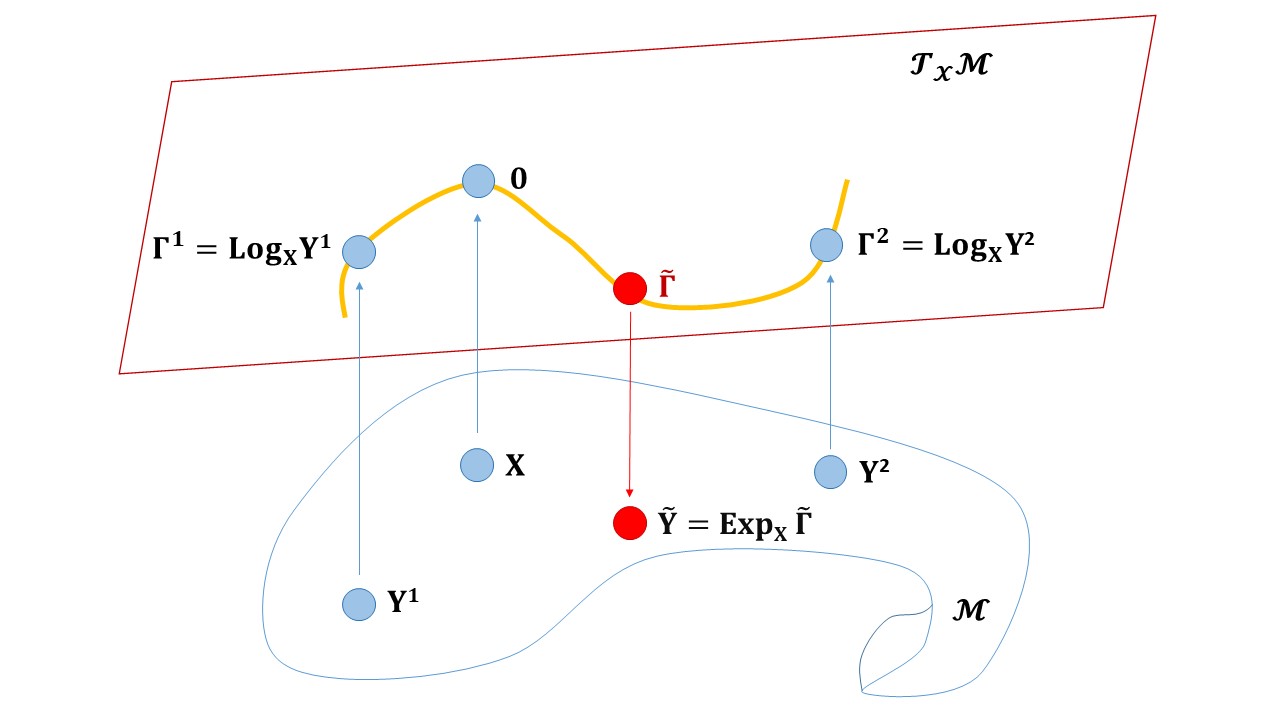}
\end{center}
\caption{Interpolation on a matrix manifold: tangent space at a reference point and logarithm and exponential maps.}
\label{fig:manifoldInterpolation}
\end{figure}

Let $\Xbold \in \mathcal{M}$ denote a reference point on the manifold $\mathcal M$, which in this case is differentiable.
Let also $\left\{\Ybold^j \in \mathcal{M}\right\}_{j=1}^{N_{\mathcal DB}}$ 
denote a set of additional points on the manifold located in a neighborhood of $\Xbold$, 
${\mathcal N}(\Xbold)$. Here and throughout the remainder of this section, the superscript $j$ refers to 
$\mubold^j \in \mathcal D \subset {\mathbb R}^{N_{\mathcal D}}$.
Then, $\Gammabold^j = \text{Log}_\Xbold(\Ybold^j)$ is a mapping of $\Ybold^j \in \mathcal{M}$ onto the 
tangent space to $\mathcal{M}$ at $\Xbold$, ${\mathcal T}_{\Xbold}\mathcal M$. The neighborhood $\mathcal{N}(\Xbold)$ is a subset of 
$\mathcal{M}$ defined such that the equation $\text{Exp}_\Xbold(\Gammabold^j) = \Ybold^j$ has a unique solution satisfying 
$\text{Log}_\Xbold(\Ybold^j) = \Gammabold^j$, where $\text{Exp}_\Xbold(\Gammabold)$ is known as the exponential map: it maps elements 
of the tangent space ${\mathcal T}_{\Xbold}\mathcal M$ back onto the manifold $\mathcal M$. In other words, the scope of 
$\mathcal{N}(\Xbold)$ is such that the inverse operation of the logarithm map has a well-defined unique output. The tangent space 
$\mathcal{T}_{\Xbold}\mathcal{M}$ being a vector space, any multi-variate interpolation approach ($N_{\mathcal D} \ge 1$)
that is valid in a vector space can be used to 
interpolate the quantities ${\Gammabold}^j$. Let $\widetilde \Gammabold$ denote the outcome of such an interpolation in 
${\mathcal T}_{\Xbold}\mathcal M$. Next, $\widetilde{\Gammabold}$ is mapped back onto $\mathcal{M}$ via the exponential map to obtain
the final outcome of the interpolation process, $\widetilde{\Ybold}(\tilde{\mubold}) = \text{Exp}_\Xbold(\widetilde{\Gammabold})$.
This process is illustrated in Figure \ref{fig:manifoldInterpolation}. 

Table \ref{ta:manifoldInterpolation} gives the expressions of the logarithm and exponential maps for various matrix manifolds 
of practical importance. The reader can observe that these expressions can be evaluated in real time.
Algorithm \ref{al:manifoldInterpolation} summarizes the procedure of interpolation on a matrix manifold
described above. Again, the reader can observe that this algorithm can be processed in real time. It can be furthermore encapsulated in the following expression
\begin{equation}\label{eq:compactInterpolation}
\widetilde{\Ybold} = \Ybold(\tilde{\mubold}) = \text{Exp}_\Xbold \left[ \mathcal{I} \left( \tilde{\mubold}; \left\{ \text{Log}_\Xbold (\Ybold^j)\right\}_{j=1}^{N_{DB}} \right) \right]
\end{equation}
where $\tilde \mubold \in \mathcal D$ denotes as before a queried but unsampled point of the parameter space $\mathcal D$,
$\mathcal{I}$ is a multi-variate interpolation operator that intervenes in the tangent space ${\mathcal T}_{\Xbold}\mathcal M$,
and $\widetilde \Ybold$ is the interpolated reduced-order matrix of interest.
(For further details on the theory and practice of interpolation on a matrix manifold in the context of linear PROMs, the reader is 
referred to \cite{amsallem2011online, amsallem2008interpolation}.)

   \begin{table}[!ht]
        \caption{Logarithm and exponential maps for some well-known matrix manifolds.} 
        \centering 
        \begin{tabular}{|c| c c c|} 
          \hline\Xhline{2\arrayrulewidth}
		\rule{0pt}{3ex}  Manifold & $\mathcal{R}^{M\times N}$ & Nonsingular matrices & SPD matrices \\ [0.25ex] 
          \Xhline{2\arrayrulewidth}
		\rule{0pt}{3ex}  $\text{Log}_{\Xbold}(\Ybold)$      &  $\Ybold-\Xbold$      &  $\text{log}(\Ybold\Xbold^{-1})$  & $\text{log}(\Xbold^{-1/2}\Ybold\Xbold^{-1/2})$ \\ [1ex]
		\Xhline{\arrayrulewidth}
		\rule{0pt}{3ex}$\text{Exp}_{\Xbold}(\Gammabold)$  &  $\Xbold+\Gammabold$  &  $\text{exp}(\Gammabold)\Xbold$   & $\Xbold^{1/2}\text{exp}(\Gammabold)\Xbold^{1/2}$ \\ [1ex]
          \Xhline{2\arrayrulewidth} 
        \end{tabular}
        \label{ta:manifoldInterpolation} 
      \end{table}

 \begin{algorithm}[!ht]
      \caption{Interpolation on a matrix manifold $\mathcal{M}$.}
      \label{al:manifoldInterpolation} 
	 \textbf{Input:} $N_\mathcal{DB}$ matrices $\Ybold^1,\ldots,\Ybold^{N_\mathcal{DB}}$ belonging to a matrix manifold $\mathcal{M}$ ; queried but unsampled parameter point $\tilde{\mubold} \in \mathcal D \subset {\mathbb R}^{N_{\mathcal D}}$\\
      \textbf{Output:} Interpolated matrix $\widetilde \Ybold \in \mathcal M$
      \begin{algorithmic}[1]
        \STATE Select a reference point on $\mathcal M$: for example, $\Xbold = \Ybold^1$
        \FOR {$j=1,\ldots,N_\mathcal{DB}$} 
           \STATE Compute $\Gammabold^j = \text{Log}_{\Xbold}(\Ybold^j)$
        \ENDFOR
        \STATE Compute each entry of $\widetilde \Gamma$ by interpolating the corresponding entries of $\Gammabold^j$, 
	      $j=1,\ldots,N_{\mathcal DB}$
	      \STATE Compute $\widetilde{\Ybold} = \Ybold(\widetilde{\mubold}) = \text{Exp}_{\Ybold^1}(\widetilde{\Gammabold})$
      \end{algorithmic}
      \end{algorithm}

\subsection{Interpolation of reduced-order model sensitivities}

The solution of the reduced-order optimization problem \eqref{eq:ROMPb} by a gradient-based algorithm requires the computation of the sensitivities of $\Abold_r^{\star}(\mubold)$ and 
$\bbold_r^{\star}(\mubold)$ with respect to the design parameters $\mu_i$, $i = 1, \cdots, N_{\mathcal D}$. Two approaches can be considered for computing these sensitivities in real time,
at a queried but unsampled parameter point $\tilde{\mubold} \in \mathcal D \subset {\mathbb R}^{N_{\mathcal D}}$:
\begin{itemize}
	\item {\it Finite differencing approach}. For each parameter component ${\tilde \mu}_i$, $i = 1, \cdots, N_{\mathcal D}$, interpolate each set of pre-computed and transformed reduced-order 
		matrices $\left\{\Abold_r^{\star}(\mubold^j)\right\}_{j=1}^{N_{\mathcal{DB}}}$ and $\left\{\bbold_r^{\star}(\mubold^j)\right\}_{j=1}^{N_{\mathcal DB}}$ on an appropriate matrix manifold 
		$\cal M$, at the parameter points ${\tilde \mubold} + {\boldsymbol \varepsilon}_i$ {\it and/or} ${\tilde \mubold} - {\boldsymbol \varepsilon}_i$, where 
		${\boldsymbol \varepsilon}_i \in {\mathbb R}^{N_{\mathcal D}}$ is the 
		vector whose entries are zero except the $i$-th entry which satisfies $\varepsilon_i \ne 0$ and $|\varepsilon_i \ll {\tilde \mu}_i|$. If the aforementioned interpolations are performed
		at the point ${\tilde \mubold} + {\boldsymbol \varepsilon}_i$ {\it or} ${\tilde \mubold} - {\boldsymbol \varepsilon}_i$ only, interpolate also each aforementioned set of 
		reduced-order matrices on $\cal M$, at the parameter point $\tilde \mubold$. Then, approximate $\displaystyle{ \frac{\partial \Abold_r}{\partial \mu_i}(\tilde\mubold)}$ and 
		$\displaystyle{ \frac{\partial \bbold_r}{\partial \mu_i}(\tilde\mubold)}$ using finite differencing.
	\item {\it Analytical approach}. For each parameter component ${\tilde \mu}_i$, $i = 1, \cdots, N_{\mathcal D}$, interpolate each set of pre-computed and transformed reduced-order 
		sensitivity matrices $\displaystyle{\left\{\frac{\partial \Abold_r^{\star}}{\partial \mu_i}(\mubold^j)\right\}_{j=1}^{N_{\mathcal DB}}}$ and 
		$\displaystyle{\left\{\frac{\partial \bbold_r^{\star}}{\partial \mu_i}(\mubold^j)\right\}_{j=1}^{N_{\mathcal DB}}}$ on an appropriate matrix manifold $\cal M$, 
		at the parameter point ${\tilde \mubold}$.
\end{itemize}
The first approach outlined above is known to lack robustness with respect to the choice of the small quantity $\varepsilon_i$. At each queried but unsampled parameter point 
$\tilde{\mubold} \in \mathcal D \subset {\mathbb R}^{N_{\mathcal D}}$, it incurs $2(1 + N_{\mathcal D})$ interpolations if first-order forward or backward differencing is adopted, 
and $2(2N_{\mathcal D})$ interpolations if second-order central differencing is used instead. The second approach mentioned above is independent of $\varepsilon_i$: it
requires $2N_{\mathcal D}$ interpolations. Hence, it is the preferred approach. However, it requires the development of an appropriate algorithm for interpolating each set
$\displaystyle{\left\{\frac{\partial \Abold_r^{\star}}{\partial \mu_i}(\mubold^j)\right\}_{j=1}^{N_{\mathcal DB}}}$ and 
$\displaystyle{\left\{ \frac{\partial \bbold_r^{\star}}{\partial \mu_i}(\mubold^j)\right\}_{j=1}^{N_{\mathcal DB}}}$ on an appropriate matrix manifold, at $\tilde \mubold$, in order to approximate
$\displaystyle{\frac{\partial \Abold_r}{\partial \mu_i}}(\tilde \mubold)$ and $\displaystyle{\frac{\partial \bbold_r}{\partial \mu_i}(\tilde \mubold)}$, respectively
\big(for the same reason explained above for the interpolation of the other components of the tuple $\mathcal{A}_r^{\star}(\mubold)$\big). Such an algoritm is proposed below, using the same notation as 
in Section \ref{sec:LAEMS}. 

First, note that each quantity $\text{Log}_\Xbold(\Ybold^j)$, where $\Ybold(\mubold^j) \in \cal M$ may represent here the $N_{w_{L_r}} \times N_{w_{L_r}}$ reduced-order matrix $\Abold_r^{\star}(\mubold^j)$ or 
the $N_{w_{L_r}} \times 1$ reduced-order matrix $\bbold_r^{\star}(\mubold^j)$ and $j = 1, \cdots, N_{\mathcal{DB}}$, depends on $\mubold^j \in \mathcal D$ but does not depend on the 
arbitrary parameter point $\mubold$. Hence, the partial derivative \eqref{eq:compactInterpolation} with respect to an arbitrary parameter component $\mu_i$ can be computed at any arbitrary
parameter point $\mubold \in \mathcal D \subset {\mathbb R}^{N_{\mathcal D}}$ as
      \begin{equation}\label{eq:interpolant_derivative}
      \frac{\partial \Ybold}{\partial \mu_i}(\mubold) 
	      = \frac{\partial \text{Exp}_{\Xbold}\mathcal{I}}{\partial \mathcal{I}}\,\frac{\partial\mathcal{I}}{\partial\mu_i}\left ({\mubold};\{\text{Log}_\Xbold(\Ybold^j)\}_{j=1}^{N_{\mathcal {DB}}}\right ), 
	      \qquad i = 1, \cdots, N_{\mathcal D} 
      \end{equation}
The objective of the algorithm presented here is to approximate the above partial derivative at a queried but unsampled parameter point $\tilde \mubold$ using interpolation on a matrix manifold. To this end,
Table \ref{ta:manifoldInterpolationDerivative} provides, for several practical and relevant matrix manifolds, the partial derivatives of the exponential map
$\displaystyle{\frac{\partial \text{Exp}_{\Xbold}\Gammabold(\mubold)}{\partial \mu_i}}$, $i = 1, \cdots, N_{\mathcal D}$. Note that 
$\displaystyle{\frac{\partial\Gammabold(\mubold)}{\partial \mu_i}} \in {\mathcal T}_{\Xbold}\mathcal M$ can be 
computed analytically, but the evaluation of the partial derivative $\displaystyle{\frac{\partial\text{exp}\left(\Gammabold(\mubold)\right)}{\partial\mu_i}}$ requires special attention. 
In order to compute this derivative, consider the matrix
\begin{equation}\label{eq:setB} 
	\Bbold_i\left(\Gammabold(\mubold)\right) = \bmat{\Gammabold(\mubold) & \displaystyle{\frac{\partial \Gammabold(\mubold)}{\partial \mu_i}} \\ \rule{0pt}{3ex} \zerobold & \Gammabold(\mubold)} \in \mathbb{R}^{2N_{w_{L_r}}\times 2N_{w_{L_r}}}, \qquad i = 1, \cdots, N_{\mathcal D} 
\end{equation}
In \cite{najfeld1995derivatives}, it is shown that the exponential of such a matrix is the matrix
\begin{equation}\label{eq:exponentialMatrixB} 
	\text{exp}\left(\Bbold_i\left(\Gammabold(\mubold)\right)\right) = \bmat{\text{exp}\left(\Gammabold(\mubold)\right) & \displaystyle{\frac{\partial \text{exp}\left(\Gammabold(\mubold)\right)}{\partial \mu_i}} \\ \rule{0pt}{3ex}\zerobold & \text{exp} \left(\Gammabold(\mubold)\right)} \in \mathbb{R}^{2N_{w_{L_r}}\times 2N_{w_{L_r}}}, \, i = 1, \cdots, N_{\mathcal D} , \, i = 1, \cdots, N_{\mathcal D} 
\end{equation}
From (\ref{eq:setB}) and (\ref{eq:exponentialMatrixB}), it follows that $\displaystyle{\frac{\partial \text{exp}\left(\Gammabold(\mubold)\right)}{\partial \mu_i}}$ can be simply obtained by 
computing the exponential of the matrix $\Bbold_i\left(\Gammabold(\mubold)\right)$, then extracting from $\text{exp}\left(\Bbold_i\left(\Gammabold(\mubold)\right)\right)$ its $(1,2)$-block. 
This computation can be performed in real time as $\Bbold_i \in {\mathbb R}^{2N_{w_{L_r}}\times 2N_{w_{L_r}}}$ is a reduced-order matrix.

\begin{table}[!ht]
\caption{Shape sensitivities of the exponential maps associated with some well-known matrix manifolds.}
\centering
\begin{tabular}{|c| c c c|}
\hline\Xhline{2\arrayrulewidth}
\rule{0pt}{3ex} Manifold & $\mathcal{R}^{M\times N}$ & Nonsingular matrices & SPD matrices \\ [0.25ex]
\Xhline{2\arrayrulewidth}
	\rule{0pt}{4ex} $\displaystyle{\frac{\partial \text{Exp}_{\Xbold}\left(\Gammabold(\mubold)\right)}{\partial \mu_i}}$  &  $\displaystyle{\frac{\partial \Gammabold(\mubold)}{\partial \mu_i}}$  &  $\displaystyle{\frac{\partial\text{exp}\left(\Gammabold(\mubold)\right)}{\partial\mu_i}\Xbold}$  & $\displaystyle{\Xbold^{1/2}\frac{\partial \text{exp}\left(\Gammabold(\mubold)\right)}{\partial\mu_i}\Xbold^{1/2}}$\\[2.0ex]
\Xhline{2\arrayrulewidth}
\end{tabular}
\label{ta:manifoldInterpolationDerivative}
\end{table}

Therefore, the proposed real-time algorithm for interpolating the quantity $\displaystyle{\frac{d \Ybold}{d \mu_i}(\mubold)}$ (\ref{eq:interpolant_derivative}) on a matrix manifold, at a queried but 
unsampled parameter point $\tilde \mubold$, where $\Ybold(\mubold) \in \cal M$ may represent the $N_{w_{L_r}} \times N_{w_{L_r}}$ reduced-order matrix $\Abold_r^{\star}(\mubold)$ or the $N_{w_{L_r}} \times 1$ 
reduced-order matrix $\bbold_r^{\star}(\mubold)$ \big(and therefore $\displaystyle{\frac{d \Ybold}{d \mu_i}}(\mubold)$ may represent the reduced-order sensitivity matrix 
$\displaystyle{\frac{d \Abold_r^{\star}}{d \mu_i}}(\mubold)$ or $\displaystyle{\frac{d \bbold_r^{\star}}{d \mu_i}}(\mubold)$\big) can be described as follows:
\begin{itemize}
	\item Let $\cal M$ be the matrix manifold to which $\Ybold(\mubold)$ belongs.  
	\item Select a reference point on $\mathcal M$: for example, $\Xbold = \Ybold^1$.
	\item Compute $\Gammabold^j = \text{Log}_{\Ybold^1}(\Ybold^j) \in {\mathcal T}_{\Ybold^1}\mathcal M$, $j = 1, \cdots, N_{\mathcal {DB}}$.
	\item Compute also analytically $\displaystyle{\frac{\partial \Gammabold(\mubold^j)}{\partial \mu_i}}$, $j = 1, \cdots, N_{\mathcal {DB}}$.
	\item Interpolate $\left\{\Gammabold^j\right\}_{j=1}^{N_{\mathcal {DB}}}$ in the tangent space ${\mathcal T}_{\Ybold^1}\mathcal M$, at the point $\tilde \mubold \in \mathcal D$,
		to obtain an approximation of $\Gammabold(\tilde \mubold)$. 
	\item For each parameter component $\tilde \mu_i$:
		\begin{itemize}
			\item Form the reduced-order matrix $\Bbold_i\left(\Gammabold(\tilde \mubold)\right)$ (\ref{eq:setB}).  
			\item Compute the reduced-order matrix $\text{exp}\left(\Bbold_i\left(\Gammabold(\tilde \mubold)\right)\right)$ (\ref{eq:exponentialMatrixB}).  
			\item Extract the $(1,2)$-block of $\text{exp}\left(\Bbold_i\left(\Gammabold(\tilde\mubold)\right)\right)$ and identify $\displaystyle{\frac{d \Ybold}{d \mu_i}(\tilde \mubold)}$ 
				with this block.
		\end{itemize}
\end{itemize}

Note that the above algorithm interpolates the values of the reduced-order sensitivity matrices $\displaystyle{\frac{\partial\Ybold}{\partial\mu_i}(\tilde \mubold)}$, $i = 1, \cdots, N_{\mathcal D}$,
on the same matrix manifold where 
Algorithm \ref{al:manifoldInterpolation} interpolates the reduced-order matrix $\Ybold(\tilde \mubold)$. Furthermore, this algorithm shares with Algorithm \ref{al:manifoldInterpolation} most
of its steps. Hence, both algorithms can be combined to interpolate on the same matrix manifold $\Ybold(\tilde \mubold)$ and its shape sensitivities 
$\displaystyle{\frac{\partial \Ybold}{\partial\mu_i}(\tilde \mubold)}$, $i = 1, \cdots, N_{\mathcal D}$.

\subsection{Tessellation-free approach based on adaptive least-squares radial basis functions}
\label{sec:RBF}

It remains to specify the interpolation method to be applied in the tangent space ${\mathcal T}_{\Ybold^1}\mathcal M$. The most important requirement for such a method is to be practical
and perform well in the presence of high-dimensional parameter spaces and scattered data. Tessellation methods or related methods based on computational meshes constitute good candidates for this
purpose when $N_{\mathcal D} \leq 3$ (for example, see \cite{amsallem2010towards}). However, they are impractical if not impossible when $N_{\mathcal D} > 3$. High-dimensional
parameter spaces call for mesh-free interpolation methods. An approach that is both popular and effective in this context is one that is based on radially symmetric basis functions that usually take the 
form of univariate functions of an Euclidean norm. Such an approach turns a multi-variate interpolation problem into one that is virtually one-dimensional, thereby simplifying the interpolation
problem and guaranteeing some interesting mathematical properties.

In the context of the parameter space $\mathcal D$ of dimension $N_{\mathcal D}$, a Radial Basis Function (RBF) is a real-valued function $\phi$ whose value at a parameter point $\mubold$ depends only 
on the distance from this point to some other point $\mathbf {c}$ called the center, so that
\begin{equation*}
\phi ( \mubold )=\phi (\Vert\mubold -\mathbf {c} \Vert_2)
\end{equation*}
Popular RBFs include the Gaussian RBF characterized by $\phi(\mubold) = \displaystyle{e^{-(\epsilon\Vert \mubold -\mathbf {c} \Vert)^2}}$, and the inverse quadratic RBF characterized by
$\phi(\mubold) = \left(1 + (\epsilon \Vert \mubold -\mathbf {c} \Vert)^2\right)^{-1}$, where $\epsilon$ is a calibration scalar. In this work, linear combinations of these RBFs are proposed
to approximate a QoI such as a component of the tuple of low-dimensional operators $\mathcal{A}_r(\tilde{\mubold})$ (\ref{eq:tuple}). The center points $\cbold$ are chosen as usual
to be the interpolation data, and therefore the interpolant can be written as
\begin{equation*}
	P(\mubold )=\sum _{j=1}^{N_{\lambdabold}}\lambda_j\,\phi (\Vert\mubold - \mubold^{j}\Vert_2), \qquad \mubold \in  \mathcal D \subset \mathbb{R}^{N_{\mathcal D}}
\end{equation*}
where the number of basis functions $N_{\lambdabold}$ is in general less than or equal to the number sampled parameter points -- that is, $N_{\lambdabold} \le N_{\mathcal{DB}}$. The coefficients 
$\lambda_j$, $j = 1, \cdots, N_{\lambdabold}$ are typically determined by solving the linear least-squares problem
\begin{equation*}
	\min_{\lambdabold \hspace{1mm}\in {\mathbb R}^{N_{\lambdabold}}} \Vert \textbf{f} - B \lambdabold \Vert_2^2
\end{equation*}
where $B \in \mathbb{R}^{N \times N_{\lambdabold}}$ is the interpolation matrix defined by $B_{ij} = \phi ( \Vert \mubold^i - \mubold^j \Vert)$, 
$\lambdabold = (\lambda_1, \ldots, \lambda_{N_{\lambdabold}})$, $\textbf{f} = (f_1, \cdots, f_N)$ is the vector storing the values of each component of the tuple $\mathcal{A}_r^{\star}({\mubold^j})$ 
to be interpolated, $j=1, \cdots, N_{\mathcal{DB}}$, and therefore $N = N_{\mathcal DB} N_{w_{L_r}} M_{w_{L_r}}$ (see Section \ref{sec:LAEMS}).

As already stated, the center points of the RBFs are typically chosen to be the entire set of sampled parameter points: this defines the straightforward implementation of an RBF-based interpolation
method. For high-dimensional parameter spaces however, or when higher accuracy is desired, this implementation may incur an unacceptable computational cost, particularly in the context
of real-time computations. Here, this cost is alleviated by focusing on the \emph{adaptive least-squares} implementation \cite{Fasshauer_RBF} of an RBF-based interpolation, where a smaller, more
economical interpolant is iteratively constructed. In this adaptive implementation, an initial interpolant is constructed using a small subset of points in $\mathcal D$ as center points. Then, the same
greedy procedure used to construct the database $\mathcal DB$ (see Section \ref{sec:database}) but equipped with the interpolation error as an error indicator, is used to select additional center points 
for constructing the RBF interpolant. Though this process incurs the construction of several small interpolants, interpolation accuracy can be achieved using a substantially smaller number of center
points than $N_{\mathcal {DB}}$, which makes this implementation more efficient than its straightforward counterpart.

\section{Applications}

In this section, the proposed fast approach for gradient-based constrained optimization based on a database of linear PROMs is 
illustrated with two examples pertaining to the aeroelastic design optimization of the same wing, under the same linearized 
flutter constraint. The only difference between the two examples is the dimension of the parameter space
$N_{\mathcal D}$ -- that is, the number of design parameters. All HDM-based computations are performed on 32 CPU cores (or simply, CPUs)
of a Linux cluster. All HDM-based computations incurred by the offline construction of a PROM database are also performed on 32 CPUs.
However, all PROM-based online computations are performed on a single CPU of the same cluster. In all cases, a CPU timing is computed and reported as the
product of the corresponding wall-clock timing and the number of CPUs used to perform the corresponding computation.

Throughout the remainder of this section, the superscript $^{\star}$ is dropped from (\ref{eq:rotatedDatabase}) and any related 
quantity for the sake of simplicity. Instead, it is assumed that all information pre-computed offline is rotated as described in
Section \ref{sec:ROT} before it is stored in a database $\mathcal {DB}$, in order to enforce mathematical consistency.

\subsection{Aeroelastic design optimization of a wing subject to a linearized flutter constraint}

Flutter is one of the most important considerations in aircraft design. It is an aeroelastic instability that can be predicted by 
linearizing the governing fluid-structure equations of dynamic equilibrium around the steady-state equilibrium associated with the 
specified flight conditions, and solving an associated eigenvalue problem. When the fluid subsystem is modeled using CFD -- for 
example, to improve fidelity in the transonic regime where shocks play an important role -- the prediction process becomes very CPU 
intensive. For this reason, the aeroelastic design optimization of an aircraft wing under a linearized flutter constraint is a good 
candidate problem for the computational technology described in this paper.

\subsubsection{Problem formulation and specification}
  
To this end, the following constrained, aeroelastic, design optimization problem -- formulated here for any generic 
aircraft wing -- is considered:
    \begin{equation}\label{eq:optimizationARW2}
    \begin{aligned}
      & \underset{\mubold\in\mathcal{D}}{\text{maximize}}  & & \frac{L(\mubold)}{D(\mubold)}  
          \\ & \text{subject to} & & W(\mubold) \leq W_{\text{upper}} 
                    \\ & & & \sigma_{\text{VM}}(\mubold) \leq \sigma_{\text{upper}}
                    \\ & & & \mubold_{\text{lower}} \leq \mubold \leq \mubold_{\text{upper}}
                    \\ & & & \zetabold_{\text{lower}}  \le \zetabold(\mubold)
    \end{aligned}
    \end{equation}
    where: $L(\mubold)$ and $D(\mubold)$ denote the lift and drag, respectively;
    $W(\mubold)$ denotes the weight of the wing, which is required to remain below an upper bound $W_{upper}$;
    $\sigma_{\text{VM}}(\mubold)$ is the von Mises stress in the wing at the steady-state equilibrium point around
    which linearization is performed, and is constrained not to exceed the maximum allowable yield stress $\sigma_{\text{upper}}$;
    $\mubold_{\text{lower}}$ and $\mubold_{\text{upper}}$ are lower and upper bounds, respectively, that define bound constraints 
    on $\mubold$ in order to avoid unrealistic designs; $\zetabold(\mubold)$ is the vector of damping ratios, where each damping ratio is required to 
    remain above a lower bound $\zeta_{\text{lower}} > 0$ in order to avoid flutter; and $\zetabold_{\text{lower}} = \zeta_{\text{lower}} {\textbf 1}$, 
    where $\textbf 1$ denotes the unit vector of same dimension as $\zetabold$ . Here, the flutter constraint is evaluated using a PROM database.
    
    Two types of design parameters are considered and therefore $\mubold$ is written as $\mubold = (\mubold_s,\mubold_m)$: aerodynamic shape parameters 
    represented by the subvector $\mubold_s\in\mathbb{R}^{N_{{\mathcal D}_s}}$, and structural parameters represented by
    $\mubold_m\in\mathbb{R}^{N_{{\mathcal D}_m}}$. The aerodynamic shape parameters $\mubold_s$ affect both the shape of the structure, 
    and the fluid flow through the transmission conditions enforced at the fluid-structure interface. The parameter subvector
    $\mubold_m$ contains the material properties of the structure as well as the internal shape parameters associated with 
    its ``dry'' elements -- that is, the elements of the structure that are not in contact with the external flow.

    The constraints of problem (\ref{eq:optimizationARW2}) can be grouped into two sets, based on how they are evaluated:
    \begin{itemize}
      \item Static aeroelastic constraints that are computed using fluid and structural
	      HDMs that are based on the three-field formulation of fluid-structure
	    interaction problems \cite{farhat1995mixed, App3F}. These constraints govern the weight, lift to drag ratio, and the von Mises stress.
            They are obtained from the solution of the governing nonlinear PDEs, which are not reduced in this paper.
      \item Dynamic, aeroelastic, flutter constraints, which are computed by linearizing the aforementioned HDMs
	    around an equilibrium state \cite{lesoinne01, lesoinne2001cfd-based}. Here, the proposed PROM database strategy is exploited
	    to alleviate the large computational cost associated with the computation of these constraints. Precisely, the computation
	    of the vector of damping ratios $\zetabold$ and its sensitivities $\displaystyle{\frac{d \zetabold}{d \mubold}}$ is accelerated using
	    the PROM interpolation approach described in Section \ref{sec:interpolation}.
    \end{itemize}
    
\begin{figure}[!ht] 
\begin{center}
\includegraphics[scale=0.110]{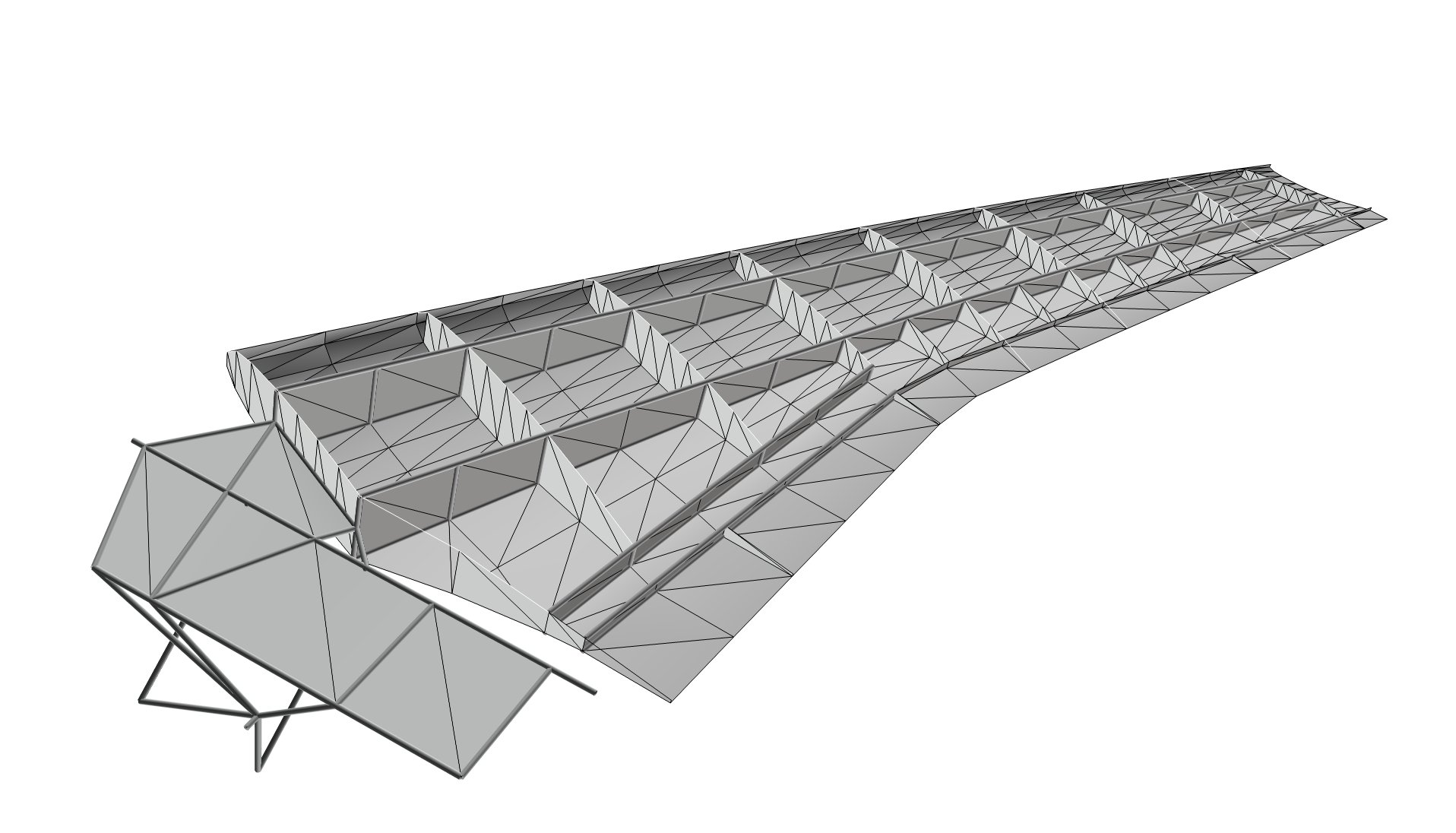}
\includegraphics[scale=0.110]{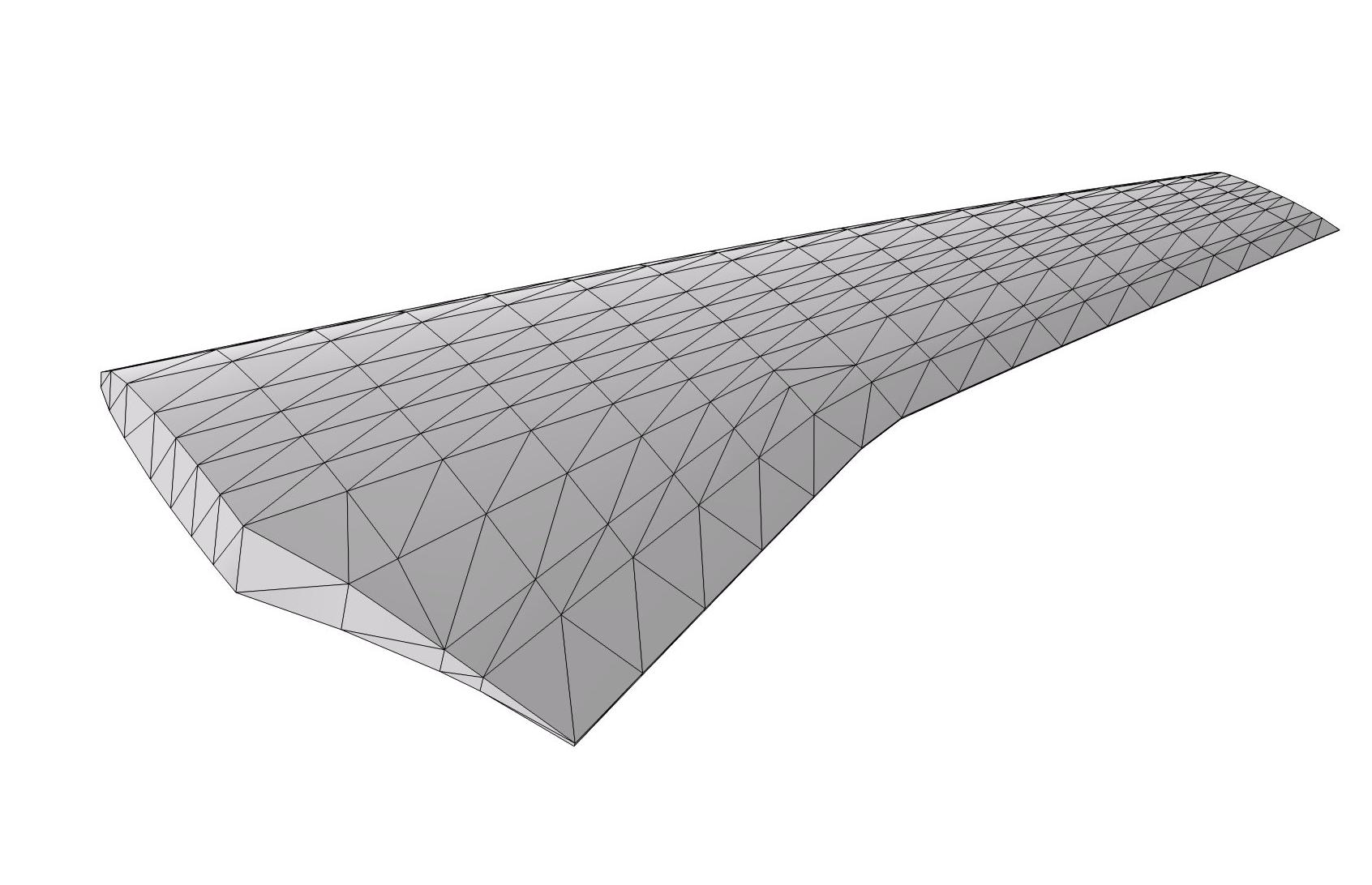}
\includegraphics[scale=0.150]{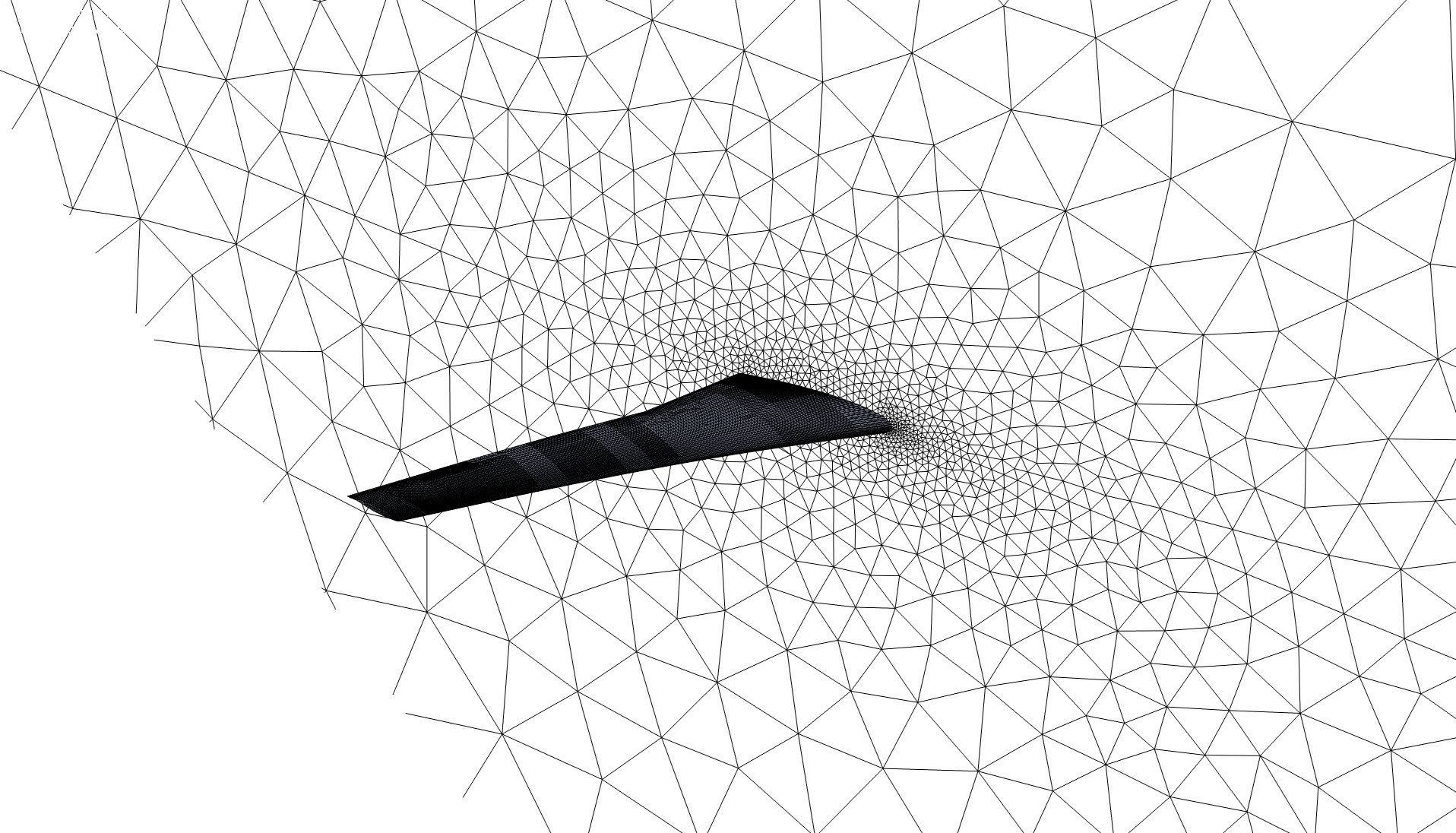}
\end{center}
\caption{Computational models for the ARW-2: detailed FE structural model with upper skin hidden (top, left); skin surface of this model (top, right); and CFD mesh (bottom).}
\label{fig:ARW2-Mesh}
\end{figure}

    Specifically, the constrained design optimization problem \eqref{eq:optimizationARW2} is considered here for a configuration
    of the Aeroelastic Research Wing (ARW-2) \cite{sandford1989geometrical}. The physical dimensions and material properties 
    of this wing, which has been built by NASA, are reported in Table \ref{ta:specARW-2}. A detailed, Finite Element (FE),
    structural model of this wing is shown in Figure \ref{fig:ARW2-Mesh}: it includes, among others, FE representations of
    the spars, ribs, hinges, and control surfaces of this wing. The model has a total of $2,556$ degrees of freedom (dofs).  
    Air is assumed to behave as a perfect gas, and its flow around this wing is assumed to be inviscid. Hence, the fluid
    flow is modeled here by the Euler equations. The computational fluid domain is discretized by a three-dimensional, unstructured, 
    CFD mesh (Figure \ref{fig:ARW2-Mesh}) with $63,484$ nodes. The focus is set on a cruise condition defined by: the altitude 
    $h = 4,000$ ft where the atmospheric density is 
    $\rho_\infty = 1.0193\times 10^{-7}\hspace{3pt} \text{lb}\cdot \text{s}^2/\text{in}^4$ and the atmospheric pressure 
    is $P_\infty = 12.7\hspace{3pt} \text{lb/in}^2$; the free-stream Mach number $M_\infty = 0.8$; and the zero angle of attack. 

\begin{table}[!ht]
        \caption{Geometrical and material properties of the ARW-2.}
        \centering
        \begin{tabular}{l c}
         \hline\Xhline{2\arrayrulewidth}
          Parameter                &  Type/Value      \\ [0.5ex]
          \hline\Xhline{2\arrayrulewidth}
          Geometry (inches)                    &                        \\ [1ex]
          \hspace{5pt}Wingspan                 &   104.9                \\ [1ex]
          \hspace{5pt}Root                     &   40.2                 \\ [1ex]
          \hspace{5pt}Tip                      &   12.5                 \\ [1ex]
          \hline\Xhline{2\arrayrulewidth}
          Material properties                  &                        \\ [1ex]
          \hspace{5pt}Skin (except flaps)      &  Composite materials   \\ [1ex]
          \hspace{5pt}Stiffeners               &  Aluminum              \\ [1ex]
          \hspace{5pt}E (psi)                  &   $1.03\times10^7$     \\ [1ex]
          \hspace{5pt}$\rho$ (lbf$\cdot s^2/\text{in}^4)$&  $2.6\times 10^{-4}$    \\ [1ex]
          \hspace{5pt}$\nu$                    &    0.32                \\ [1ex]
          \Xhline{2\arrayrulewidth}
        \end{tabular}
      \label{ta:specARW-2}
\end{table}

For this wing, the upper bound for the von Mises stress constraint, $\sigmabold_{\text{upper}}$, is set to $2.5\times10^4$ psi, and 
that for the weight, $W_{\text{upper}}$, is set to $400$ lbs. Several values of the lower bound for the damping ratios, 
$\zeta_{\text{lower}}$, are considered in the interval $[4.1\times 10^{-4}, 4.5\times 10^{-4}]$ in order to study how the optimal 
aeroelastic design varies with this flutter threshold.

\begin{figure}[!ht]
\begin{center}
\includegraphics[scale=0.25]{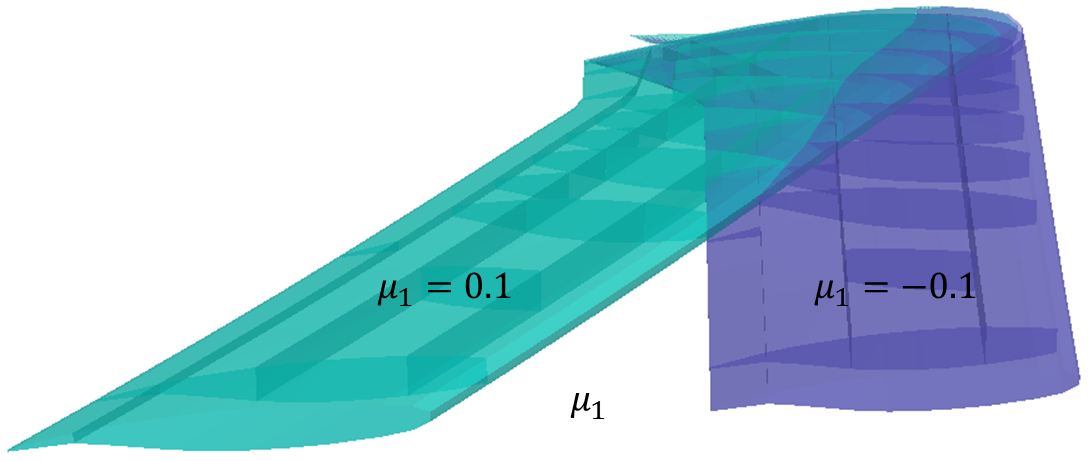} 
\includegraphics[scale=0.25]{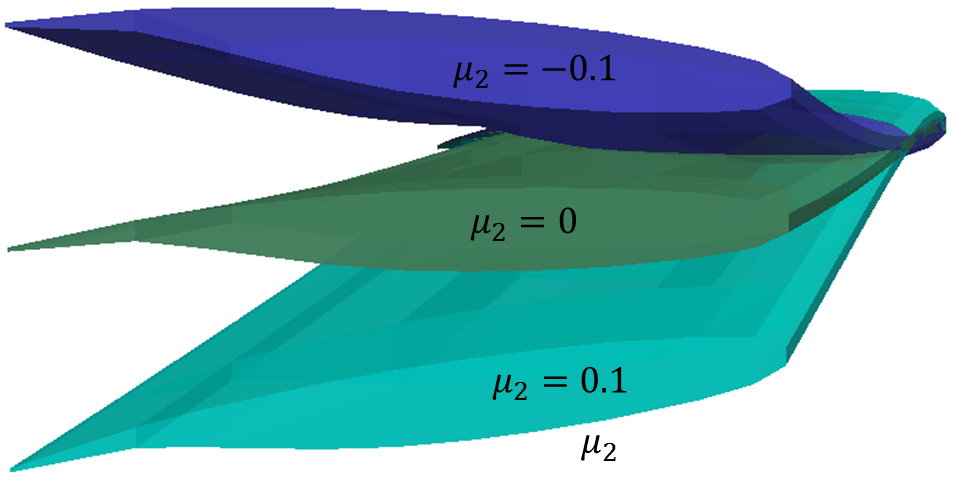} 
\includegraphics[scale=0.25]{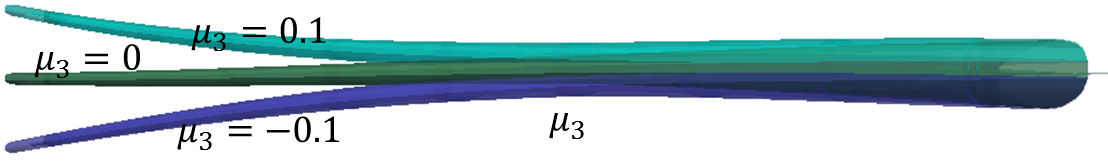} 
\end{center}
\caption{Aerodynamic shape design parameters.} 
\label{fig:ARW2parametervariations}
\end{figure}

Three shape and three structural parameters are considered -- that is, $\mubold_s \in \mathbb{R}^3$
and $\mubold_m\in\mathbb{R}^3$. The shape parameters are the back sweep angle $\mu_s^1$, the twist angle $\mu_s^2$, and the 
dihedral angle $\mu_s^3$. These parameters, which are graphically depicted in Figure \ref{fig:ARW2parametervariations},
are constrained to have a magnitude less than 0.1 radian. 

\begin{figure}[!ht]
\begin{center}
\includegraphics[scale=0.2]{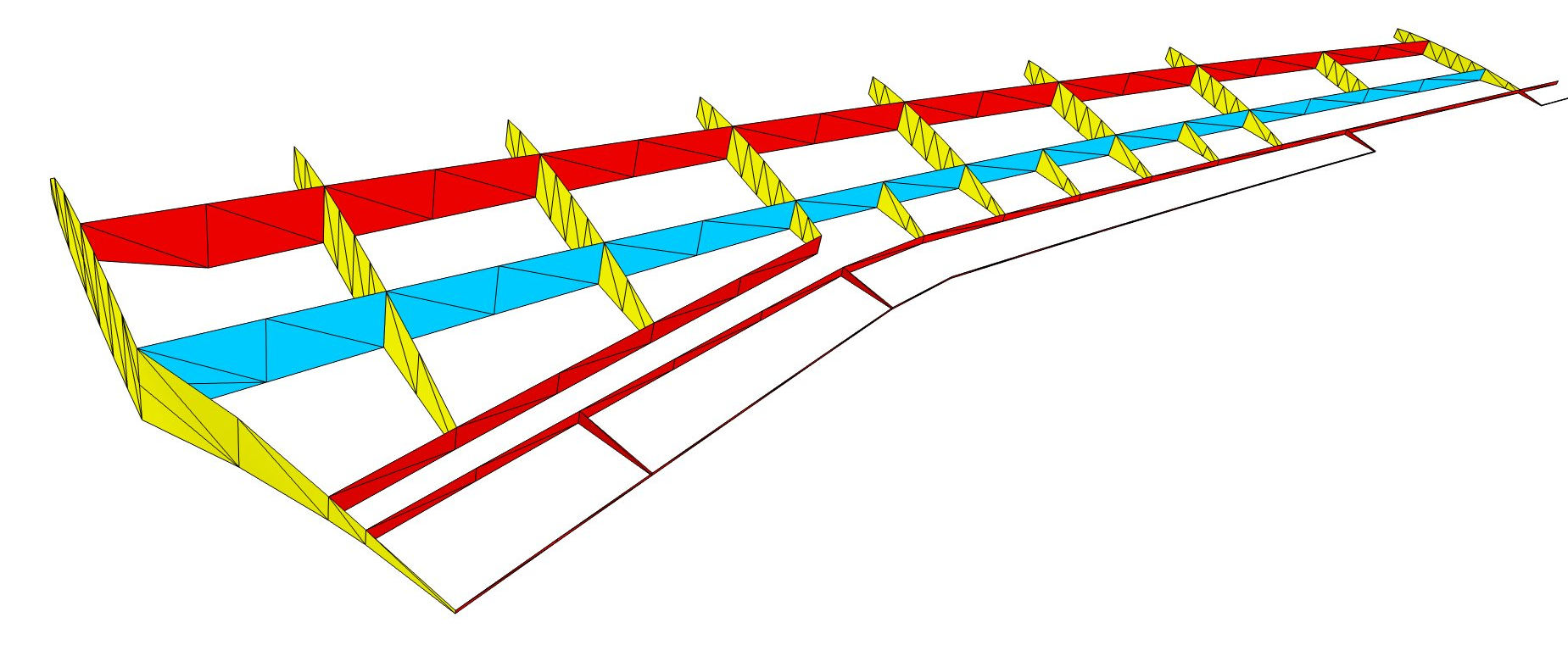} 
\end{center}
\caption{Organization of the ARW-2 stiffeners in three different groups.}
\label{fig:ARW2stiffenersGroup}
\end{figure}

For the purpose of optimization, the stiffeners of the ARW-2 are organized in three different groups:
the outer spars, the central spar, and the ribs (see Figure \ref{fig:ARW2stiffenersGroup}). Each of these groups contains
elements with different thicknesses, but is attributed a single, collective thickness increment as an optimization parameter.
Hence, the structural parameters stored in $\mubold_m$ are the three thickness increments for the outer spars ($\mu_m^1$), 
the central spar ($\mu_m^2$), and the ribs ($\mu_m^3$). The magnitude of each of these increments is constrained to be less 
than 0.01 in -- which corresponds on average to 10\% of the thickness of a stiffening element of the considered initial configuration
of the ARW-2.

\vspace{10 mm}
      
\vspace{10 mm}

Specifically, two different instances of the design optimization problem (\ref{eq:optimizationARW2}) are considered in the remainder of this section for the ARW-2: 
\begin{enumerate}
\item A first instance where only the aerodynamic shape parameters are activated as optimization parameters -- that is, $\mubold = \mubold_s$.  
\item A second instance where both the aerodynamic shape and structural parameters are activated as optimization parameters -- that is, $\mubold = (\mubold_s,\mubold_m)$.
\end{enumerate}
For each instance outlined above, the performance of the model reduction approach adopted here for computing the flutter constraints and its impact on the overall performance of the NAND solution algorithm is measured relative 
to the performance results obtained using the HDM-based solution strategies.
  
In order to keep this paper as self-contained as possible, the three-field formulation of fluid-structure interaction problems and its linearization, which are used here to formulate the flutter constraint 
$\zetabold(\mubold) \geq \zetabold_{\text{lower}}$ in (\ref{eq:optimizationARW2}), are next overviewed (for more details, the reader is referred to \cite{farhat1995mixed, App3F, lesoinne01, lesoinne2001cfd-based}.
  
\subsubsection{Linearized CFD-based computational aeroelasticity and model order reduction}\label{sec:linearizedFSI}
  
A high-fidelity, nonlinear, CFD-based fluid-structure (or aeroelastic) system can be described by a three-field Arbitrary Lagrangian Eulerian (ALE) formulation \cite{farhat1995mixed, App3F} that treats the moving CFD mesh as a
pseudo-structural system \cite{farhat98}. After semi-discretization, linearization of the governing semi-discrete equations about 
a fluid-structure equilibrium state \cite{lesoinne01, lesoinne2001cfd-based}, and elimination from these equations of the CFD mesh 
position vector, this formulation leads to the following system of linear Ordinary Differential Equations (ODEs)
\begin{eqnarray}
\Cbold_v(\mubold) {\dot \wbold} + \Hbold(\mubold)  {\wbold} + \Rbold(\mubold) {\dot\ubold} + \Gbold(\mubold) {\ubold} &=& \zerobold\label{eq:FsysLinParam}\\
\Mbold(\mubold) \ddot\ubold + \Dbold(\mubold) \dot\ubold + \Kbold(\mubold) \ubold &=& \Pbold(\mubold) \wbold\label{eq:SsysLinParam}
\end{eqnarray}
where: a dot denotes a derivative with respect to time $t$; $\Cbold_v\in\mathbb{R}^{N_f\times N_f}$ is the diagonal matrix of 
cell volumes in the CFD mesh, and $N_f$ is the dimension of the semi-discretized fluid subsystem;  $\wbold(t)\in\mathbb{R}^{N_f}$
denotes the conservative fluid state vector; $\ubold(t)\in\mathbb{R}^{N_s}$ is the vector of structural displacements of 
dimension $N_s$; $\Hbold \in\mathbb{R}^{N_f \times N_f}$ and $\Gbold \in\mathbb{R}^{N_f \times N_x}$ denote the Jacobians of the 
numerical fluid fluxes with respect to $\wbold$ and $\xbold$, respectively; $\Rbold \in\mathbb{R}^{N_f \times N_x}$ is the Jacobian of the
numerical fluid fluxes with respect to the fluid mesh velocity; $\Mbold\in\mathbb{R}^{N_s\times N_s}$ is the mass matrix associated 
with the FE representation of the structural subsystem; $\Dbold\in\mathbb{R}^{N_s\times N_s}$ and 
$\Kbold\in\mathbb{R}^{N_s\times N_s}$ are the FE damping and stiffness matrices, respectively; and 
$\Pbold\in\mathbb{R}^{N_s\times N_f}$ denotes the Jacobian of the aerodynamic forces with respect to $\wbold$ and acts
on the ``wet'' surface of the structure -- that is, the fluid-structure interface.

First, the dimensionality of the structural subsystem is reduced using modal truncation. For this purpose,
a ROB $\Xbold(\mubold)\in\mathbb{R}^{N_s\times k_s}$, where $\Xbold$ depends on the parameter vector
$\mubold = (\mubold_s,\mubold_m)\in\mathbb{R}^{p_s+p_m}$, is constructed using the first $k_s$ modes of the structural subsystem.
Then, the subspace approximation
\begin{equation*}\label{eq:approxS}
\ubold(t) \approx \Xbold(\mubold) \ubold_r(t)
\end{equation*}
where $\ubold_r\in\mathbb{R}^{k_s}$ denotes the vector of $k_s$ generalized coordinates associated with the modal approximation,
is performed. Finally, Equation \eqref{eq:SsysLinParam} is reduced by Galerkin projection onto $\Xbold(\mubold)$ to obtain
\begin{equation*}\label{eq:SROM}
\ddot\ubold_r + \Dbold_r(\mubold) \dot \ubold_r + \Omegabold_r^2(\mubold) \ubold_r 
= \Xbold(\mubold)^T\Pbold(\mubold) \wbold
\end{equation*}
where $\Dbold_r(\mubold) = \Xbold(\mubold)^T\Dbold(\mubold)\Xbold(\mubold)\in\mathbb{R}^{k_s\times k_s}$ and $\Omegabold_r^2\in\mathbb{R}^{k_s\times k_s}$ is the diagonal matrix of eigenvalues associated with the $k_s$ eigenmodes.

The dimensionality of the fluid subsystem is reduced using the POD method in the frequency domain \cite{kim98}. Specifically, a 
ROB $\Vbold(\mubold)\in\mathbb{R}^{N_f\times k_f}$ is constructed using Algorithm~\ref{alg:PODfreq}. Note that by construction,
this global ROB satisfies the orthogonality condition $\Vbold(\mubold)^T\Cbold_v(\mubold)\Vbold(\mubold) = \Ibold_{k_f}$.
Next, the state vector $\wbold$ is approximated as 
\begin{equation*}\label{eq:approxF}
\wbold(t) \approx \Vbold(\mubold) \wbold_r(t)
\end{equation*}
where $\wbold_r$ denotes the vector of $k_f$ generalized coordinates associated with $\Vbold(\mubold)$, and Galerkin projection is
performed to transform (\ref{eq:FsysLinParam}) into the set of reduced-order coupled linear ODEs
\begin{eqnarray}
\dot\wbold_r + \Hbold_r(\mubold) \wbold_r 
+ \Rbold_r(\mubold) \dot\ubold_r 
+ \Gbold_r(\mubold) \ubold_r  &=& \boldsymbol{0}\label{eq:eqRedF}
\\ \ddot\ubold_r + \Dbold_r(\mubold) \dot \ubold_r 
+ \Omegabold_r^2(\mubold) \ubold_r - \Pbold_r(\mubold)\wbold_r &=& \boldsymbol{0}\label{eq:eqRedS}
\end{eqnarray}
where $\Hbold_r(\mubold) =\Vbold(\mubold)^T\Hbold(\mubold_s)\Vbold(\mubold)$, 
$\Rbold_r(\mubold) = \Vbold(\mubold)^T\Rbold(\mubold_s)\Xbold(\mubold)$, and 
$\Pbold_r(\mubold) = \Xbold(\mubold)^T\Pbold \Vbold(\mubold)$. If needed, the stability of this resulting
fluid-structure PROM can be enforced using the fast algorithm described in \cite{amsallem2014stability}.

\begin{algorithm}[t]
      \caption{Construction of a fluid ROB for an aeroelastic system.}
      \label{alg:PODfreq}
        \textbf{Input:} Parameter vector $\mubold=(\mubold_s,\mubold_m)\in\mathbb{R}^{p_s+p_m}$, 
                        frequency sampling $\{\xi_j\}_{j=1}^{N_\xi}$, 
                        structural eigenmodes $\Xbold(\mubold)\in\mathbb{R}^{N_s\times k_s}$, 
                        desired ROB dimension $k_f$ \\
        \textbf{Output:} ROB $\Vbold(\mubold)$
        \begin{algorithmic}[1] 
          \STATE Compute $\Cbold_v = \Cbold_v(\mubold)$, 
                         $\Hbold = \Hbold(\mubold)$, 
                         $\Rbold = \Rbold(\mubold)$, and 
                         $\Gbold = \Gbold(\mubold)$
          \FOR {$i=1,\cdots,k_s$}
            \FOR {$l=1,\cdots,N_\xi$}
              \STATE Solve the linear system $(j\xi_l \Cbold_v + \Hbold)\wbold_{i,l} = -(j\xi_l\Rbold+\Gbold)\xbold_i$, 
                     where $\xbold_i$ is the $i$-th vector in $\Xbold$
            \ENDFOR
          \ENDFOR
          \STATE Construct the complex-valued snapshot matrix $\Wbold = [\wbold_{1,1},\cdots,\wbold_{k_s,N_\omega}]$
          \STATE Compute the SVD $\Cbold_v^{\frac{1}{2}}[\text{Re}(\Wbold),\text{Im}(\Wbold)] = \Ubold \Sigmabold\Zbold^T$
          \STATE Retain the first $k_f$ left singular vectors $\Vbold(\mubold) = \Cbold_v^{-\frac{1}{2}}\Ubold_{k_f}$
        \end{algorithmic}
      \end{algorithm}

The coupled PROM (\ref{eq:eqRedF})--(\ref{eq:eqRedS}) can be re-written in compact form as
\begin{equation*}\label{eq:NODE1}
\dot\qbold_r = \Nbold_r(\mubold) \qbold_r
\end{equation*}
where 
\begin{equation}
\label{eq:NODE2}
\qbold_r(t) = \left[\begin{array}{c}\wbold_r(t) \\\dot\ubold_r(t) \\\ubold_r(t)\end{array}\right] \in\mathbb{R}^{k_f+2k_s}
\end{equation}
and
\begin{equation}\label{eq:Nop}
\Nbold_r(\mubold) 
= \left[\begin{array}{ccc}-\Hbold_r(\mubold) & -\Rbold_r(\mubold) & - \Gbold_r(\mubold) 
\\ \Pbold_r(\mubold) & -\Dbold_r(\mubold) & -\Omegabold_r^2(\mubold) 
\\ \boldsymbol{0} & \Ibold_{k_s} & \boldsymbol{0}\end{array}\right]\in\mathbb{R}^{(k_f+2k_s)\times(k_f+2k_s)}
\end{equation}

For simplicity, and without any loss of generality, structural damping is neglected throughout the remainder of this section
-- that is $\Dbold_r(\mubold) = \zerobold$. 

Given that $k_s \ll N_s$ and $k_f \ll N_f$, the reduced-order eigenvalue problem
\begin{equation}
\label{eq:eigenvalueproblem}
	\Nbold_r(\mubold){\qbold}_{r}(\mubold) = \lambda(\mubold)\qbold_{r}(\mubold)
\end{equation}
can be solved in real-time, for any specified value of $\mubold$, to obtain the $2k_s+k_f$ eigenpairs
$\{\lambda_j, ~ \qbold_{r_j}\}_{j=1}^{2k_s+k_f}$. For each computed eigenvalue $\lambda_j$, which can be complex-valued, the corresponding damping ratio $\zeta_{j}$ is defined as 
\begin{equation}\label{eq:dampRatioDef}
\zeta_{j} =-\frac{\lambda_j^R}{\sqrt{\left(\lambda_j^R\right)^2+\left(\lambda_j^I\right)^2}}
\end{equation}
where $\lambda_j^R$ and $\lambda_j^I$ are the real and imaginary part of $\lambda_j$, respectively.

During the solution of the optimization problem \eqref{eq:optimizationARW2}, many values of the parameter vector $\mubold$
will be queried. For any value $\tilde \mubold$ for which the blocks of the matrix $\Nbold_r(\mubold)$ (\ref{eq:Nop}) are not available
in the PROM database ${\mathcal DB}$ constructed offline, the PROM interpolation method described in 
Section \ref{sec:interpolation} is applied to construct $\Nbold_r(\tilde\mubold)$ and its sensitivities with respect to the
optimization parameters block-by-block. For this purpose, Equation (\ref{eq:eigenvalueproblem}) can be identified with
Equation (\ref{eq:PROM}) -- or more specifically, the reduced-order matrix $\Nbold_r$ can be identified with the reduced-order matrix ${\mathbf A}_r$ -- even if
(\ref{eq:LinearForm}) is a linear system of equations but (\ref{eq:eigenvalueproblem}) is a linear eigenvalue problem.

The solution of the optimization problem \eqref{eq:optimizationARW2} also requires the sensitivity of the vector
of damping ratios with respect to the optimization parameters, $\displaystyle{\frac{\partial \zetabold}{\partial \mu_i}}$,
$i = 1,~\cdots,~N_{\mathcal D}$. Appendix \ref{sec:appendix} describes a real-time algorithm for computing online these
sensitivities.

\subsection{Offline database construction}

Three different instances of the greedy procedure are applied offline to construct three different databases of coupled PROMs of the form (\ref{eq:Nop}) and of reduced dimensions $k_s = 6$ and
$k_f = 100$ and compare their performances. In all three instances: the set of candidate parameter points $\Xi$ is constructed using FF design and five points along each dimension of the parameter 
space, which leads to $N_{\Xi} = 125$; the convergence tolerance is set to $\epsilon_{tol} = 5\times 10^{-2}$; and the error indicator is chosen to be the residual-based error indicator 
(\ref{eq:residual}) with ${\mathbf R}_L$ set -- for this purpose, and only for this purpose -- to the residual associated with the first eigenpair of the PROM matrix (\ref{eq:Nop}).

The three aforementioned instances differ as follows:
\begin{itemize} 
\item \textit{Standard.} In this case the residual-based error indicator is evaluated at each iteration at every candidate parameter point in $\Xi$.
\item \textit{Random.} In this case, the residual-based error indicator is evaluated at each iteration at every candidate parameter point in a random subset $\Pi_m \subset \Xi$ of size $N_{\Pi} = 20$.
\item \textit{Saturation.} In this case, the greedy procedure is that described in Algorithm \ref{sec:saturationgreedy} with $N_{\Pi} = 20$ and $N_{\Pi}^{\prime} = 50$. The saturation constant is set to 
	$\tau_s=2$ in order to capture a fluctuation of the error at some parameter points as more of these points are sampled to refine the database of PROMs.  
\end{itemize}

The performances of these databases are contrasted in Table~\ref{ta:comparisonGreedyAlgorithms} for $\mubold=\mubold_s$. The reader can observe that the accelerated parameter sampling procedure based on 
a saturation assumption described in Section~\ref{sec:EGP} and Algorithm \ref{sec:saturationgreedy} reduces the computational cost of the standard greedy procedure by a factor greater than 7.

All performance results reported in the remainder of this paper are for the database of coupled PROMs constructed using Algorithm \ref{sec:saturationgreedy}.

\begin{table}[t]
        \caption{Performance comparison of three different instances of the greedy procedure.}
        \centering
        \begin{tabular}{l c c c}
         \hline\Xhline{2\arrayrulewidth}
            Greedy procedure& \textit{Standard} & \textit{Random} & \textit{Saturation}  \\ [0.5ex]
          \Xhline{2\arrayrulewidth}
            $N_{\Xi}$ &    125            &     125         &      125              \\ [1ex]
            $N_{\Pi}$ &    -              &      20         &       20              \\ [1ex]
          Number of sampled HDMs &   16   &      15         &       15              \\ [1ex]
          Maximum Relative Indicated Error ($\%$) &   4.9  &     3.9         &      4.6              \\ [1ex]
          Total wall-clock time (hrs) &60.9 &    12.1         &      8.4              \\ [1ex]
          Speed-up with respect to standard greedy procedure& 1&     5.1         &      7.3              \\ [1ex]
          \Xhline{2\arrayrulewidth}
        \end{tabular}
      \label{ta:comparisonGreedyAlgorithms}
      \end{table}

\subsection{Speed-up factors for one online evaluation of the flutter constraint} 
                
Table \ref{ta:singleDynamicConstrintComputation} compares the performance of the proposed methodology for accelerating the evaluation of a discretized, parametric, linear PDE constraint
with that of the HDM-based evaluation of this constraint, for a single computation of the vector of damping ratios in problem (\ref{eq:optimizationARW2}). For this purpose, this computation 
($\zetabold$ and $\displaystyle{\frac{d\zetabold}{d\mubold}}$) is carried out using the fixed-point iteration algorithm described in \cite{amsallem2015linearized}. 
For the case where $\mubold = \mubold_s$, the PROM database approach is shown to deliver a speed-up factor roughly equal to $17$ for the wall-clock time, and $535$ for the CPU time. 
For the case $\mubold = (\mubold_s,\mubold_m)$, the speed-up factor for the wall-clock time is roughly equal to $29$, and that for the CPU time is slightly larger than $926$.

\begin{table}[t]
\caption{Performance results for the evaluation of the flutter constraint and its sensitivities for a single queried parameter vector.}
\centering
\begin{tabular}{l r r r|r r r}
\hline\Xhline{2\arrayrulewidth}
Computational model & \multicolumn{3}{c|}{$\mubold=\mubold_s$ ($N_{\mathcal D}=3$)}    & \multicolumn{3}{c}{$\mubold=(\mubold_s,\mubold_m)$ ($N_{\mathcal D}=6$)} \\
\cline{2-7}
& $\#$ of CPUs & Wall-clock & CPU  & $\#$ of CPUs & Wall-clock & CPU  \\
&              & time       & time &              & time       & time \\ [0.5ex]
\Xhline{2\arrayrulewidth}
\hspace{5pt}HDM           &       32           &    0.342         & 10.934          &       32           &    0.492         & 15.744  \\ [1ex]
\hspace{5pt}Database of PROMs&        1           &    0.019         & 0.019           &        1           &    0.017         & 0.017  \\ [1ex]
\hspace{5pt}Speed-up factor &                    &    18.0          & 575.5           &                    &    28.9          & 926.1  \\ [1ex]
\Xhline{2\arrayrulewidth}
\end{tabular}
\label{ta:singleDynamicConstrintComputation}
\end{table}
       
\subsection{Speed-up factors for the solution of the optimization problem with $N_{\mathcal D} = 3$}\label{sec:3param}

\begin{figure}[H]
\begin{center}
\includegraphics[scale=0.35]{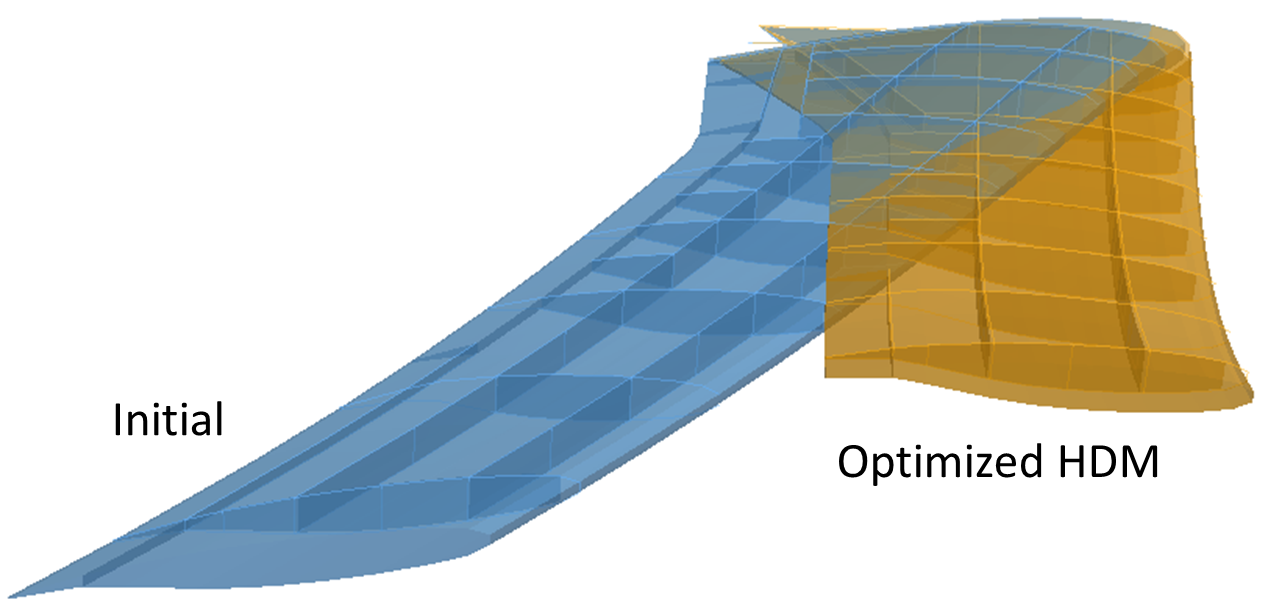} 
\end{center}
\caption{Initial and optimized designs of the ARW-2.} 
\label{fig:initOptConfigs}
\end{figure}

\begin{figure}[H]
\begin{center}
\includegraphics[scale=0.45]{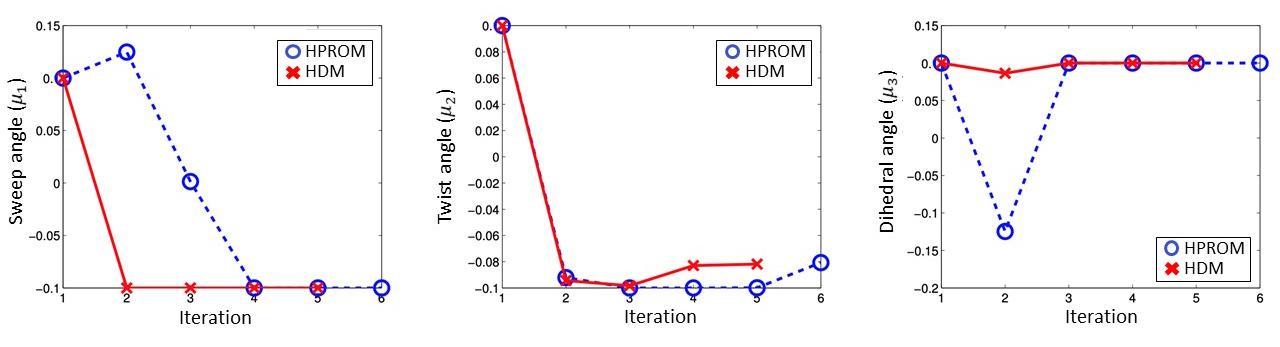}
\end{center}
\vspace{-4mm}
\caption{Iteration-histories of the shape parameters.}
\label{fig:designVariableHistoryOptimization}
\end{figure}

\begin{figure}[H]
\begin{center}
\includegraphics[scale=0.45]{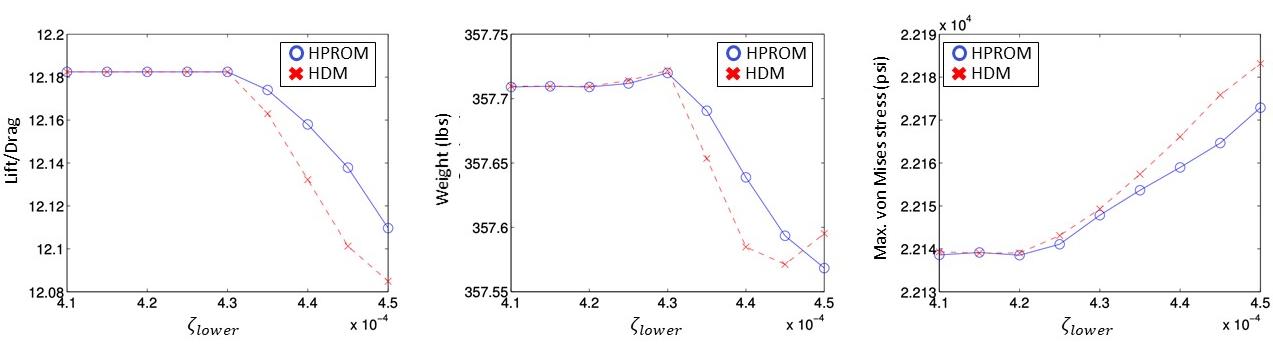}
\end{center}
\vspace{-5mm}
\caption{Iteration-histories of various QoIs.} 
\label{fig:historyOptimization}
\end{figure}

Table \ref{ta:singleOptResults} summarizes the results obtained for the solution of the optimization problem (\ref{eq:optimizationARW2}) with $\mubold=\mubold_s$
and $\zeta_{\text{lower}}=4.2\times 10^{-4}$. Given the same initial shape defined by $\mubold = (0.1,0.1,0.1)$, both computational approaches based on the HDM and
PROM database are found to deliver almost the same optimal design, after almost the same number of iterations. (The initial and optimized shapes of the ARW-2 are shown in Figure \ref{fig:initOptConfigs}).
However, Figure \ref{fig:designVariableHistoryOptimization} shows that the optimization paths followed by both computational models are different. 
In both cases, the flutter constraint is violated initially, but satisfied at convergence to the optimal shape. After optimization, the lift to drag ratio is increased by $15.1 \%$,
and the weight and maximum von Mises stresses are increased by $4.4 \%$ and $18.2 \%$, respectively, without violating their respective constraints.
Figure \ref{fig:historyOptimization} shows the iteration-histories of the lift to drag ratio, minimum damping ratio, maximum von Mises stress, and weight. Again,
this figure shows that the two computational models contrasted here converge to the same results, but follow different paths in general.

Table \ref{ta:singleOptResults} also reports the CPU timing results obtained for each considered computational model. The CPU timing result reported for the offline phase of model reduction corresponds to the CPU time 
elapsed in constructing the PROM database using the accelerated parameter sampling procedure described in Section \ref{sec:EGP}. The ``static'' constraints refer to all of the weight, von Mises, and box constraints, 
as well as the objective function. The flutter constraint refers as before to the constraint on the vector of damping ratios. The following observations are noteworthy:
\begin{itemize}
\item The CPU-based speed-up factor achieved for the online evaluation of the flutter constraint is consistent with that reported in Table \ref{ta:singleDynamicConstrintComputation} for a single such evaluation.
\item The smaller speed-up factor achieved for the entire online phase is due to the fact that the so-called static constraints are not reduced.  
\item Because it requires an offline investment, the PROM database approach can be expected to deliver a (significant) speed-up factor only if this investment is amortized. As this table shows, such an investment
is rarely amortized by the solution of a single problem -- in this case, an optimization problem. However, it will be shown next that such an investment can be amortized if multiple optimizations are performed 
-- for example, as in multi-start global optimization strategies. 
\end{itemize}
      
\begin{table}[th]
\caption{Performance results for the solution of the ARW-2 shape optimization problem with $\zeta_{\text{lower}} = 4.2\times 10^{-4}$ and $\mubold=\mubold_s$ ($N_{\mathcal D} = 3$).}
\centering
\begin{tabular}{l r r r}
\hline\Xhline{2\arrayrulewidth}
Design                              &      Initial       &  Optimized         & Optimized            \\ [0.5ex]
\Xhline{2\arrayrulewidth}
Computational model                 &                    &        HDM         & PROM database        \\ [0.5ex]
\Xhline{2\arrayrulewidth}
	\# of CPUs                  &                    &        32          &     variable         \\ [0.5ex]
\Xhline{2\arrayrulewidth}
\hspace{5pt}$\mubold_s$             & (0.1,0.1,0.1)      & (-0.1,-0.0814,0.1) & (-0.1,-0.0807,0.1) \\ [1ex]
\hspace{5pt}$L/D$                   & 10.6               &     12.2           &    12.2         \\ [1ex]
\hspace{5pt}Weight (lbs)            & 342.6              &   357.7            &    357.7        \\ [1ex]
\hspace{5pt}Min. $\zetabold$        &$2.24\times 10^{-4}$&$4.22\times 10^{-4}$&$4.25\times 10^{-4}$ \\ [1ex]
\hspace{5pt}Max. $\sigma_{\text{VM}}$ (psi)  & 18,731.9  &   22,139.5         &    22,139.0        \\ [1ex]
\hspace{5pt}\# of iterations        &                    &      5             &       6            \\ [1ex]
\hspace{5pt}\# of function evaluations &                 &     24             &      22            \\ [1ex]
\Xhline{2\arrayrulewidth}
CPU time (hours)                     &                    &                    &                  \\ [1ex]
\hspace{5pt}Offline phase (32 CPUs) &                    &       0            &   268.8          \\ [1ex]
\hspace{5pt}Online evaluation of the static constraints  (32 CPUs) &        &   14.2             &    14.1          \\ [1ex]
\hspace{5pt}Online evaluation of the flutter constraint (1 CPU)  &       &   240.5            &     0.43           \\ [1ex]
\hspace{5pt}Total (online)          &                    &   254.7            &    14.5          \\ [1ex]
\hspace{5pt}Total (offline + online)  &                    &   254.7            &   283.3          \\ [1ex]
\Xhline{2\arrayrulewidth}
Speed-up factors based on CPU time             &                    &                    &                  \\ [1ex]
\hspace{5pt}Speed-up factor for flutter constraint &                &    1               &  559.3           \\ [1ex]
\hspace{5pt}Speed-up factor for online phase       &                    &    1               &   17.6           \\ [1ex]
\hspace{5pt}Total speed-up factor              &                    &    1               &   0.7            \\ [1ex]
\Xhline{2\arrayrulewidth}
\end{tabular}
\label{ta:singleOptResults}
\end{table}

Table \ref{ta:robustOptCompTime} reports the CPU timing results obtained for multiple instances of the flutter optimization problem \eqref{eq:optimizationARW2}, corresponding to the following
nine different values of the lower bound for the damping ratios, $\zeta_{\text{lower}} = \{4.1, 4.15, 4.2, 4.25, 4.3, 4.35, 4.4, 4.45, 4.5\}\times 10^{-4}$, for the purpose of studying the influence of this flutter 
threshold on the optimized design. Furthermore for each damping ratio lower bound, a multi-start strategy with ten different initial points is applied in order to provide a more robust search for a global optimizer. Therefore, 
the performance results reported in Table \ref{ta:robustOptCompTime} correspond to 90 independent optimization problems that are solved using both of the HDM and PROM database computational models.
In this case, the reader can observe that:
\begin{itemize}
\item The proposed PROM database approach delivers a speed-up factor based on the CPU time that is equal to $1,459.3$ for the flutter constraint, and a speed-up factor based on the CPU time that is equal to $43.6$ for the overall 
online phase.
\item This approach achieves a speed-up factor based on the CPU time that is equal to $32.4$ for the entire computations, including those incurred by the construction offline of the database of PROMs, which confirms the observation made above.
\end{itemize}

\begin{figure}[H]
         \begin{center}
	    \includegraphics[scale=0.45]{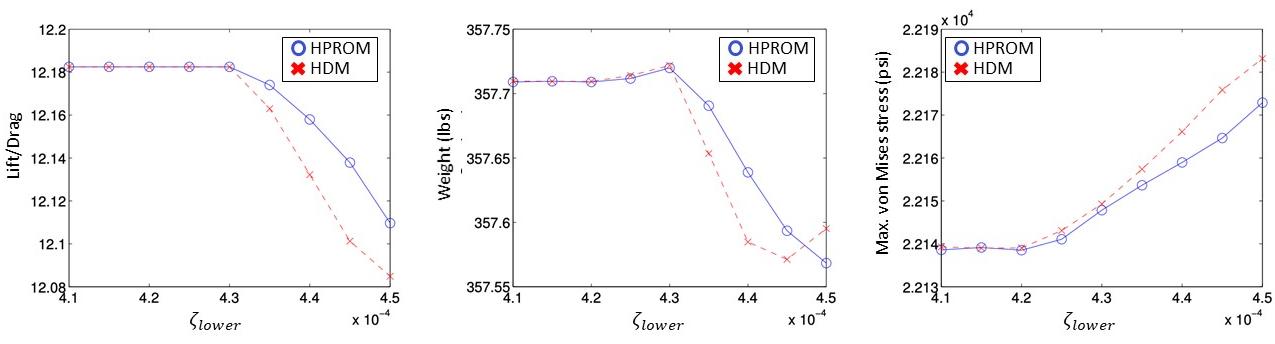}
	 \end{center}
	    \vspace{-7mm}
         \caption{Influence of the lower bound of the flutter constraint on various QoIs associated with the computed optimal design.}
         \label{fig:robustOptLDWeightStress}
       \end{figure}

\begin{figure}[H]
         \begin{center}
	    \includegraphics[scale=0.45]{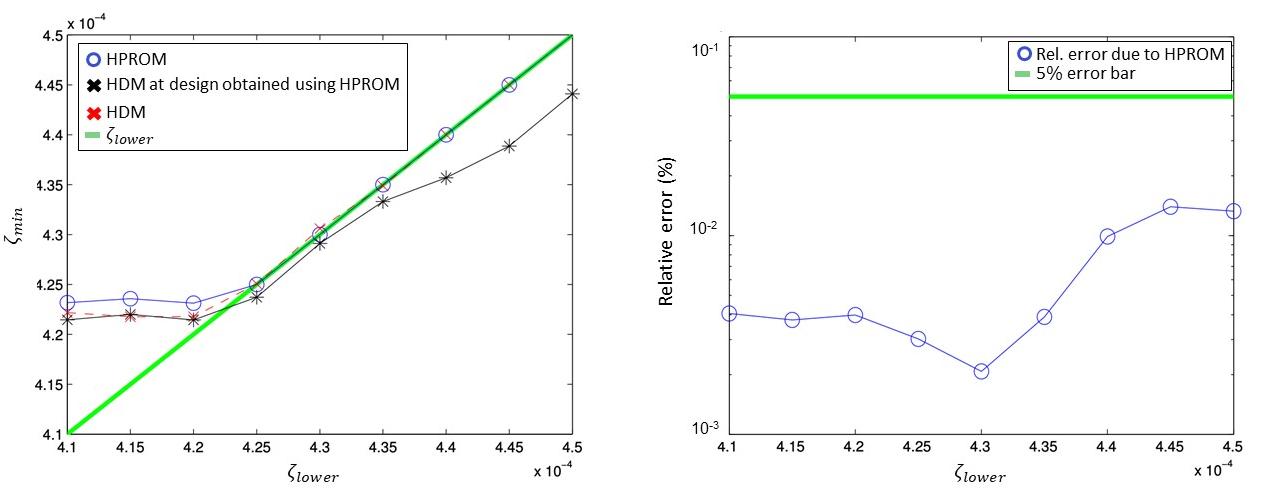}
	 \end{center}
	    \vspace{-7mm}
         \caption{Influence of the lower bound of the flutter constraint on: the minimum damping ratio of the computed optimal design (left); and the relative error associated with the computation of this QoI using the PROM database approach.}
         \label{fig:robustOptDamping}
       \end{figure}

Figure \ref{fig:robustOptLDWeightStress} graphically depicts the influence of $\zeta_{\text{lower}}$ on the optimal lift to drag ratio, the optimal weight, and the maximum von Mises stress associated in each case with the optimal
design. Overall, this influence appears to be minimal. Similarly, Figure \ref{fig:robustOptDamping}-left shows the influence of $\zeta_{\text{lower}}$ on the minimum damping ratio associated with an optimal design. 

Figure \ref{fig:robustOptDamping}-right reports on the relative error associated with the computation of the minimum damping ratio using the PROM database approach, for all considered values of $\zeta{\text{lower}}$. It shows that 
in all cases, this error is significantly below the convergence threshold of $5\%$ adopted during the construction offline of the PROM database.

      \begin{table}[H]
        \caption{Performance results for the solution of nine instances of the ARW-2 shape optimization problem with $\mubold=\mubold_s$ ($N_{\mathcal D} = 3$), but different values of $\zeta_{\text{lower}}$.}
        \centering
        \begin{tabular}{l r r}
         \hline\Xhline{2\arrayrulewidth}
          Computational model                    &        HDM    &       PROM Database   \\ [0.5ex]
	  \Xhline{2\arrayrulewidth}
	  \# of CPUs                             &        32     &     variable         \\ [0.5ex]
          \Xhline{2\arrayrulewidth}
          CPU time (hours)                       &               &                      \\ [1ex]
          \hspace{5pt}Offline phase (32 CPUs)    &      0        &   268.8              \\ [1ex]
          \hspace{5pt}Online evaluation of the static constraints (32 CPUs)         & 2,065.1       &   755.7              \\ [1ex]
          \hspace{5pt}Online evaluation of the flutter constraints (1 CPU)       &31,813.5       &    21.8              \\ [1ex]
          \hspace{5pt}Online miscellaneous       &    14.2       &    0.06              \\ [1ex]
          \hspace{5pt}Total (online)             &33,892.8       &   777.6              \\ [1ex]
          \hspace{5pt}Total (offline + online)   &33,892.8       & 1,046.4              \\ [1ex]
          \Xhline{2\arrayrulewidth}
          Speed-up factors based on CPU time     &               &                      \\ [1ex]
          \hspace{5pt}Speed-up factor for flutter constraints&      1        & 1,459.3              \\ [1ex]
          \hspace{5pt}Speed-up factor for online phase  &      1        &    43.6              \\ [1ex]
          \hspace{5pt}Total speed-up factor             &      1        &    32.4              \\ [1ex]
          \Xhline{2\arrayrulewidth}
       \end{tabular}
       \label{ta:robustOptCompTime}
       \end{table}

\subsection{Speed-up factors for the solution of the optimization problem with $N_{\mathcal D} = 6$}\label{sec:6param}

Table \ref{ta:singleOptResultsSixParams} presents the performance results obtained for the solution of the optimization problem \eqref{eq:optimizationARW2} with $\mubold=(\mubold_s,\mubold_m)$ and 
$\zeta_{\text{lower}}=4.2\times 10^{-4}$. As in the case of three shape parameters ($\mubold = \mubold_s$), both computational models based on the HDM and PROM database approaches
are found to deliver almost the same optimal shape, given the same initial design. However, the two different computational models lead in this case to two slightly different optimal structural 
parameters. In both cases, the flutter constraint is violated initially but satisfied at convergence, and the lift to drag ratio is increased by $8.0 \%$. However, the weight is increased by $4.7 \%$ in 
the optimal design found using the PROM database as a computational model, but is maintained almost constant in the counterpart design found by the HDM-based optimization process. The maximum von Mises 
stress is decreased by $6.6 \%$ in the optimal design found using the PROM database computational model, and by $6.4 \%$ in the counterpart design found using the HDM-based approach. Although the 
optimal shape parameters obtained in this case  -- $\mubold_s=(-0.1,-0.1,0.1)$ -- are similar to those obtained in the case of optimization using only the three shape parameters (see 
Table \ref{ta:singleOptResults}), the maximum von Mises stress experienced in this case by the optimal design is smaller than its counterpart experienced by the optimal design obtained when 
$\mubold=\mubold_s$. This lower value is achieved by increasing the thickness of the stiffeners without violating the weight constraint. The reader can also observe that for both considered 
computational models, similar iteration counts and similar numbers of function evaluations are performed.

Table \ref{ta:singleOptResultsSixParams} also reports the CPU time results and associated speed-up factors achieved for this instance of the optimization problem \eqref{eq:optimizationARW2}. Due to the 
larger computational expense induced by the computation of a larger number of sensitivities, the CPU time-based speed-up factor achieved by the PROM database approach for a single evaluation of the 
flutter constraint and its sensitivities is substantially higher (speed-up factor = 1,154.7) than for the case where $N_{\mathcal D}=3$ (speed-up factor = 559.3). The speed-up factor for the overall 
online phase of model reduction is similar however in magnitude. Again, the total speed-up factor for the combined offline and online phases is less than 1, due to the computational expense incurred
by building the database. As demonstrated however for the case where $N_{\mathcal D}=3$, this computational overhead can be amortized and significant overall speed-up factors can be achieved when many 
optimization problems are solved using the same PROM database. 
      
      
      \begin{table}[th]
      \caption{Performance results for the solution of the ARW-2 shape optimization problem with $\zeta_{\text{lower}} = 4.2\times 10^{-4}$ and $\mubold=(\mubold_s, \mubold_m)$ ($N_{\mathcal D} = 6$).}
        \centering
        \begin{tabular}{l r r r}
         \hline\Xhline{2\arrayrulewidth}
                Design             &      Initial                    &  Optimized         & Optimized      \\ [0.5ex]
          \Xhline{2\arrayrulewidth}
          Computational model      &                                 &        HDM         &       PROM Database  \\ [0.5ex]
          \Xhline{2\arrayrulewidth}
	  \# of CPUs               &                                 &         32         &       variable       \\ [0.5ex]
	  \Xhline{2\arrayrulewidth}
          \hspace{5pt}$\mubold_s$                & (0,0,0)           & (-0.1,-0.1,0.1)   & (-0.1,-0.1,0.1)  \\ [1ex]
          \hspace{5pt}$\mubold_m$                & (0,0,0)           & (0.1,0.1,-0.015)  & (0.1,0.1,0.013)  \\ [1ex]
          \hspace{5pt}$L/D$                      & 11.3              &    12.2           & 12.2             \\ [1ex]
          \hspace{5pt}Weight (lbs)               & 349.9             &    349.8          & 366.3            \\ [1ex]
          \hspace{5pt}Min. $\zetabold$           &$3.7\times 10^{-4}$&$4.2\times 10^{-4}$&$4.2\times 10^{-4}$\\ [1ex]
          \hspace{5pt}Max. $\sigma_{\text{VM}}$ (psi)& 20,297.1      &  18,991.9         & 18,953.3         \\ [1ex]
          \hspace{5pt}\# of iterations           &                   &      6            &     5                 \\ [1ex]
          \hspace{5pt}\# of function evaluations &                   &      11           &     9                 \\ [1ex]
          \Xhline{2\arrayrulewidth}
          CPU time (hours)                       &                   &                   &                       \\ [1ex]
          \hspace{5pt}Offline phase (32 CPUs)    &                   &      0            &   8,800               \\ [1ex]
          \hspace{5pt}Online evaluation of the static constraints (32 CPUs)&          &     11.0          &   12.8                \\ [1ex]
          \hspace{5pt}Online evaluation of the flutter constraint (1 CPU) &          &     173.2         &   0.15                \\ [1ex]
          \hspace{5pt}Total (online)             &                   &     184.2         &   13.0                \\ [1ex]
          \hspace{5pt}Total (offline + online)   &                   &     184.2         &   8,813               \\ [1ex]
          \Xhline{2\arrayrulewidth}
          Speed-up factors based on CPU time     &                   &                   &                       \\ [1ex]
          \hspace{5pt}Speed-up factor for flutter constraint&                   &      1            &   1,154.7              \\ [1ex]
          \hspace{5pt}Speed-up factor for online phase &                   &      1            &   14.2                \\ [1ex]
          \hspace{5pt}Total speed-up factor      &                   &      1            &   0.02                \\ [1ex]
          \Xhline{2\arrayrulewidth}
       \end{tabular}
       \label{ta:singleOptResultsSixParams}
       \end{table}

\section{Conclusions}
\label{sec:Conclusions}
A novel methodology is introduced in this paper for accelerating the solution of Partial Differential Equation (PDE)-constrained optimization problems where at least one PDE-based constraint is linear. 
This methodology consists in constructing offline a database of pointwise, linear, Projection-based Reduced-Order Models (PROMs), and equipping it with two essential computational technologies: an
offline pre-processing algorithm for enforcing mathematical consistency among the pre-computed PROMs that are stored in the database; and an online algorithm for interpolating on matrix manifolds the 
pre-computed PROMs and their sensitivities at any queried but unsampled parameter point in the design parameter space. To reduce the computational overhead incurred by the construction offline of the 
PROM database, a parameter sampling procedure based on an appropriate saturation assumption is proposed to minimize the number of pointwise PROMs to be pre-computed and stored in the database, while 
maximizing the accuracy of the online interpolation of these PROMs at any queried but unsampled parameter point. The interpolation on matrix manifolds, which preserves some desirable properties
of the PROMs to be interpolated and their sensitivities, incurs the standard interpolation in the tangent space to a considered manifold -- which is a vector space -- of some parametric 
quantities. It is shown in this paper that when the design parameter space is high-dimensional, such an interpolation can be practically and effectively performed using radial basis functions.
A new and rigorous approach for interpolating on matrix manifolds matrix sensitivities is also presented. The PROM database approach and its peripherals are illustrated with the solution of
a number of different instances of a PDE-constrained, aeroelastic, optimization problem associated with NASA's ARW-2 wing and involving a linear, PDE-based flutter constraint. Using this
problem as a background problem, it is also shown in this paper that when multiple optimizations are performed -- for example, as in multi-start global optimization strategies -- the computational 
overhead incurred by the construction offline of a PROM database can be amortized and a significant overall, CPU time-based, speed-up factor can be achieved. Finally, it is noted that given the
recent work published in \cite{washabaugh2016deform}, the methodology presented in this paper is readily extendible to the case where the PDE-based optimization problem contains a combination
of linear and nonlinear PDE-based constraints.

 \section{Appendix}\label{sec:appendix} 
 \paragraph{Real-time computation of the sensitivity of the vector of damping ratios with respect to the optimization parameters, $\displaystyle{\frac{\partial \zetabold}{\partial \mu_i}}$,
 $i = 1,~\cdots,~N_{\mathcal D}$}
For each eigenpair $(\lambda_j, \qbold_{r_j})$, $j = 1,~\ldots,~k_f+2k_s$, the differentiation of \eqref{eq:eigenvalueproblem} with respect to $\mu_i,~i=1,\cdots,N_{\mathcal D}$, leads to
      \begin{equation*}\label{eq:derivativeofreduced}
        \frac{\partial \Nbold_r}{\partial \mu_i} \qbold_{r_j} - 
         \frac{\partial \lambda_j}{\partial \mu_i} \qbold_{r_j} + 
        (\Nbold_r - \lambda_j \Ibold)\frac{\partial \qbold_{r_j}  }{\partial \mu_i} = \zerobold
      \end{equation*}
where the reduced-order matrix matrix $\Nbold_r$ and the associated reduced-order vector $\qbold_r$ are defined in (\ref{eq:Nop}) and (\ref{eq:NODE2}), respectively, $\qbold_{r_j}$ is the $j$-th
right eigenvector of $\Nbold_r$, $\lambda_j$ is the corresponding eigenvalue, $\Ibold$ is the identity matrix, and $\mu_i$ is the $i$-th component of the parameter vector $\mubold \in \mathcal D \subset \Rbb^{N_{\mathcal D}}$.
Multiplying (\ref{eq:derivativeofreduced}) from the left by $\pbold_{r_j}^H$, where $\pbold_{r_j}$ denotes the $j$-th left eigenvector of $\Nbold_r$ and the superscript $H$ designates the transpose of the conjugate of a complex-valued
quantity, and noting that $\pbold_{r_j}^T(\Nbold_r - \lambda_j \Ibold)=\zerobold$ gives
      \begin{equation*}\label{eq:structuraleigenvaluesensitivity}
        \frac{\partial \lambda_j}{\partial \mu_i} = 
       \frac{\displaystyle{\pbold_{r_j}^H\frac{\partial \Nbold_r}{\partial \mu_i}\qbold_{r_j}}}{\displaystyle{\pbold_{r_j}^H \qbold_{r_j}}}
      \end{equation*} 
In the above result, the sensitivity matrix $\displaystyle{\frac{\partial \Nbold_r}{\partial \mu_i}}$ can be computed at any queried but unsampled parameter point in real-time, block-by-block, using interpolation on matrix manifolds
as described in Section \ref{sec:linearizedFSI} (see also Table \ref{ta:manifoldInterpolation} and Table \ref{ta:manifoldInterpolationDerivative}). 
Then, from Equation (\ref{eq:dampRatioDef}) and the chain rule, it follows that the sensitivity of each damping ratio $\zeta_j$ with respect to any parameter component $\mu_i$ can be computed in real-time as follows
      \begin{equation}\label{eq:dampingratioderivative}
        \frac{\partial\zeta_{j}}{\partial \mu_i} = \frac{\partial\lambda_j^R}{\partial \mu_i}   
                    \left(\frac{(\lambda_j^R)^2}{|\lambda_j|^3} - \frac{1}{|\lambda_j|}\right) + 
        \frac{\partial\lambda_j^I}{\partial \mu_i}  \frac{\lambda_j^R \lambda_j^I}{|\lambda_j|^3}, \quad j = 1,~\ldots,~k_f+2k_s, \quad i = 1,~\ldots,~N_{\mathcal D}
      \end{equation}
      
\section{Acknowledgements}
Youngsoo Choi acknowledges partial support by the office of Naval Research under Grant N00014-11-1-0707 and partial support by the Army High Performance Computing Research Center under 
Cooperative Agreement W911NF-07-2-0027, while in residence at Stanford University. Gabriele Boncoraglio and Charbel Farhat acknowledge partial support by the Office of Naval Research under 
Grant N00014-17-1-2749, partial support by the Boeing Company under Contract Sponsor Ref. 134824, and partial support by a research grant from the King Abdulaziz City for Science and Technology (KACST). 
Spenser Anderson acknowledges the support by the Department of Defense through the National Defense Science \& Engineering Graduate Fellowship (NDSEG) Program.
This document does not necessarily reflect the position of these institutions, and no official endorsement should be inferred.

\section{References}

\bibliographystyle{model1-num-names}
\bibliography{pap}
\end{document}